\documentstyle[twoside,11pt,leqno]{article}

 \def\X{{\cal X}} \def\S{{\cal
S}} \def\H{{\cal H}} 
\def\D{{\texttt{D}}}
\def\i{\rm{iso}} 
\def\n{{\texttt N}} 
\def\C{\texttt{C}} \def\iso{\textrm{iso}}

\def\B{B({\cal H})} \def\b{B({\cal X})}
\def\asc{ \textrm{asc}} \def\dsc{ \textrm{dsc}}

\newtheorem{df}{Definition}[section]
\newtheorem{thm}[df]{Theorem} \newtheorem{pro}[df]{Proposition}
\newtheorem{cor}[df]{Corollary} 
\newtheorem{rema}[df] {Remark} \newtheorem{lem}[df] {Lemma}
\def\sfstp{{\hskip-1em}{\bf.}{\hskip1em}}

\def\subject#1{\renewcommand{\thefootnote}{}\footnote
{AMS(MOS) subject classification (2010). Primary: {#1}}}

\def\keywords#1{\renewcommand{\thefootnote}{}\footnote
{Keywords: {#1}}}

\def\enddemo{\qed \endtrivlist} \expandafter\let\csname
enddemo*\endcsname=\enddemo

\def\qedsymbol{\ifmmode\bgroup\else$\bgroup\aftergroup$\fi
\vcenter{\hrule\hbox{\vrule
height.5em\kern.5em\vrule}\hrule}\egroup}
\def\qed{\ifmmode\else\unskip\nobreak\fi\quad\qedsymbol}

\title{\bf  Structure of $n$-quasi left $m$-invertible and related classes of operators}
\author{\normalsize B.P. Duggal, I.H. Kim}
\date{}

\begin{document}

\maketitle \thispagestyle{empty} \vskip-16pt

\subject{Primary47A05, 47A55 Secondary47A80, 47A10.}
\footnote{The second named author was supported by Basic Science Research Program through the National Research Foundation of Korea (NRF)
 funded by the Ministry of Education (NRF-2019R1F1A1057574).}
\begin{abstract}  Given Hilbert space operators $T, S\in\B$, let $\triangle$ and $\delta\in B(\B)$ denote the elementary operators $\triangle_{T,S}(X)=(L_TR_S-I)(X)=TXS-X$ and $\delta_{T,S}(X)=(L_T-R_S)(X)=TX-XS$. Let $d=\triangle$ or $\delta$.  Assuming $T$ commutes with $S^*$, and choosing $X$ to be the positive operator $S^{*n}S^n$ for some positive integer $n$, this paper exploits properties of elementary operators to study the structure of $n$-quasi $[m,d]$-operators  $d^m_{T,S}(X)=0$ to bring together, and improve upon, extant results for a number of classes of operators, amongst them $n$-quasi left $m$-invertible operators, $n$-quasi $m$-isometric operators, $n$-quasi $m$-selfadjoint operators and $n$-quasi $(m,C)$ symmetric operators (for some conjugation $C$ of $\H$).  It is proved that  $S^n$ is the perturbation by a nilpotent of the direct sum of an operator $S_1^n=(S|_{\overline{S^n(\H)}})^n$ satisfying $d^m_{T_1,S_1}(I_1)=0$, $T_1=T|_{\overline{S^n(\H)}}$, with the $0$ operator; if also $S$ is left invertible,  then $S^n$ is similar to an operator $B$ such that $d^m_{B^*,B}(I)=0$.  For power bounded $S$ and $T$ such that  $ST^*-T^*S=0$ and $\triangle_{T,S}(S^{*n}S^n)=0$,  $S$ is  polaroid (i.e., isolated points of the spectrum are poles). The product property, and the perturbation by a commuting nilpotent property, of operators $T, S$ satisfying $d^m_{T,S}(I)=0$, given certain commutativity properties, transfers to operators satisfying $S^{*n}d^m_{T,S}(I)S^n=0$. 
\end{abstract}

%%%%%%%%%%%%%%%%%%%%%%%%%%%%%%%%%%%%%%%% SECTION 1

\section {\sfstp Introduction} Let $\b$ (resp., $\B$) denote the algebra of operators, equivalently bounded linear transformations, on a complex infinite dimensional Banach space $\X$ (resp., Hilbert space $\H$) into itself. Given operators $T, S\in\b$, let $L_T$ and $R_S\in B(\b)$ denote, respectively,  the operators $$L_T(X)=TX,  R_S(X)=XS$$ of left multiplication by $T$ and right multiplication by $S$. The elementary operators $\triangle_{T,S}$ and $\delta_{T,S}\in B(\b)$ are then defined by  $$\triangle_{T,S}(X)=(L_TR_S-I)(X)=TXS-X$$ and  $$\delta_{T,S}(X)=(L_T-R_S)(X)=TX-XS.$$  Let $d_{T,S}\in B(\b)$ denote either of the operators $\triangle_{T,S}$ and $\delta_{T,S}$. Let $I$ denote the identity of $\b$ and let $m\geq 1$ be some integer. Then  
 \begin{eqnarray} \label{eqnarray1}\triangle^m_{T,S}(I)=\triangle_{T,S}(\triangle^{m-1}_{T,S})(I))=\sum_{j=0}^m{(-1)^{j}\left(\begin{array}{clcr}m\\j\end{array}\right)T^{m-j}S^{m-j}} \end{eqnarray}
and  \begin{eqnarray} \label{eqnarray1}\delta^m_{T,S}(I)=\delta_{{T,S}}(\delta^{m-1}_{T,S}(I))=\sum_{j=0}^m{(-1)^{j}\left(\begin{array}{clcr}m\\j\end{array}\right)T^{m-j}S^j} \end{eqnarray}
We say in the following that an operator $S\in\b$ is an $m-(d,T)$ operator if $d^m_{T,S}(I)=0$. Examples of $m-(d,T)$  operators $S\in\b$ occur quite naturally. Thus:  if an operator $S\in\b$ is {\it $m$-left invertible} by $T\in\b$, then $$\triangle_{T,S}^m(I)=\sum_{j=0}^m{(-1)^{j}\left(\begin{array}{clcr}m\\j\end{array}\right)T^{m-j}S^{m-j}}=0$$ \cite{{OAMS},{DM},{CG}}; if $S\in\b$ is {\it $m$-isometric}, then $$\triangle_{S^*,S}^m(I)=\sum_{j=0}^m{(-1)^{j}\left(\begin{array}{clcr}m\\j\end{array}\right) S^{*(m-j)}S^{m-j}}=0$$ \cite{{AS},{BMMN},{D}}; if $S\in\B$ is {\it $m$-selfadjoint}, then  $$\delta^m_{S^*,S}(I)=\sum_{j=0}^m{(-1)^{j}\left(\begin{array}{clcr}m\\j\end{array}\right) S^{*(m-j)}S^j}=0$$ \cite{L} and if $S\in\B$ is {\it $(m,C)$-isometric} for some conjugation $C$ of $\H$ , then  $$\delta^m_{S^*,CSC}(I)=\sum_{j=0}^m{(-1)^{j}\left(\begin{array}{clcr}m\\j\end{array}\right) S^{*(m-j)}CS^jC}=0$$ \cite{CKL}. Operators $S\in  m-(d,T)$, in particular the classes consisting of $m$-isometric and $(m,C)$-isometric operators \cite{CLM}, have been studied in a number of papers in the recent past (see cited references for further references). A gereralization of the class consisting of of $m$-isometric  (resp., $(m,C)$-isometric)  operators which has drawn some attention in the recent past is that of the $n$-quasi  $m$-isometric  (resp., $n$-quasi $(m,C)$-isometric) operators, where an operator $S\in\B$ is said to be $n$-quasi $m$-isometric (resp., $n$-quasi $(m,C)$-isometric) for some integer $n\geq 1$ if $S^{*n}\triangle^m_{S^*,S}(I)S^n=\triangle^m_{S^*,S}(S^{*n}S^n)=0$ (respectively, $S^{*n}\triangle_{S^*,CSC}(I)S^n=0$) \cite{ {MP},{ACL}}. In keeping with current terminology \cite{{K}, {MP}, {ACL}}, we say in the following  that an operator {\it $S\in\B$ is $n$-quasi $[m,d]$-intertwined by $T\in\B$ (equivalently, $T$ is an $n$-quasi $[m,d]$-intertwining of $S$)} for some integer $n\geq 1$ if 
$$S^{*n} d^m_{T,S}(I)S^n=0.$$  It is immediate from the definition that if $S\in\B$ is 
$n$-quasi $[m,d]$-intertwined by $T$, $[S,T^*]=ST^*-T^*S=0$ (thus $S^{*n} d^m_{T,S}(I)S^n=d^m_{T,S}(S^{*n}S^n)=0$), $T^*_1=T^*|_{\overline{S^n(\H)}}$
and $S_1=S|_{\overline{S^n(\H)}}$, then $d^m_{T_1,S_1}(I_1)=0$. Choosing $T=S^*$, we prove in the following that if $S^{*n}d^m_{S^*,S}(I)S^n=0$ and  if $d=\triangle$ (resp., $d=\delta$ and $S$ is injective), then there exist a positive operator $Q$ and an operator $A$ such that $\triangle^m_{A^*,A}(Q)=0$  and $S^n$ is similar to $A$ (resp., $\delta^m_{A^*,A}(Q)=0$ and $\delta_{A,S^n}(P)=0$, $P$ a quasi-affinity). Furthermore, if $S$ is left invertible, then there exists an operator $B\in\B$ such that $S^n$ is similar to $B$ and $d^m_{B^*,B}(I)=0$.

\vskip3pt  

 Left $m$-invertible Banach space (as also $m$-isometric, $m$-selfadjoint Hilbert space \cite{L}) operators are known to satisfy the properties that:
if $S_i, T_i\in \b$, $i=1,2$, are such that $S_i$ is left $m_i$-invertible by $T_i$ and $[S_1,S_2]=
0=[T_1,T_2]$, then $S_1S_2$ is left $(m_1+m_2-1)$-invertible by $T_1T_2$; if $S_1\in\b$ is left $m_1$-invertible 
by $T_1\in\b$ and $N_1\in\b$ is an $n_1$-nilpotent operator which commutes with $S_1$, then $S_1+N_1$ is left $(m_1+n_1-1)$-invertible by $T_1$ \cite{CG}. These results, which hold equally well for $[m,d]$-intertwinings, have extensions to $n$-quasi $[m,d]$-intertwining (Hilbert space) operators $S,T$. Let us say that {\em $S_1\in\B$ 
	is $n(S)$-quasi  $[m,d]-$intertwined by $T_1\in\B$ for some operator $S\in\B$ if 
	$$S^{*n}d^m_{T_1,S_1}(I)S^n=0.$$} We prove that if $S_i, T_i\in\B$ ($i=1,2$) are some operators such that $S_1$ is $n(S)$-quasi $[m_1,d]$-intertwined by $T_1$, $S_2$ is $[m_2,d]$-intertwined by $T_2$, $[S_1,S_2]=0=[T_1,T_2]$ and $[S,S_i]=0=[S,T_i^*]$ ($i=1,2$), then $S_1S_2$ is $n(S)$-quasi $[(m_1+m_2-1),d]$-intertwined by $T_1T_2$. For an $n$-quasi $m_1$-isometric $S\in\B$ and an $m_2$-isometric $T\in\B$ such that $S, T$ commute, this implies that $ST$ is an $n$-quasi $(m_1+m_2-1)$-isometry.  Again, if $S$ is $n(S)$-quasi $[m,d]$-intertwined by $T$, $N_i\in\B$ are nilpotent operators ($i=1,2$), $[S,N_1]=0=[S,T^*]$,  $[N_2,T]=0=[S,N_2^*]$ and $S$ is injective in the case in which $d=\delta$, then $(S^*+N^*_1)^{n+n_1-1} d^{m+n_1+n_2-2}_{T+N_2,S+N_1}(I)(S+N_1)^{n+n_1-1}=0$.  Translated to left invertible  $n$-quasi $m$-isometric operators $S\in\B$ such that $S$ commutes with an $n_1$-nilpotent operator $N\in\B$ this implies that there exists an $m$-isometric operator $B\in\B$ such $(S+N)^{n+n_1-1}$ is similar to $B$.

\

Recall that a Banach space operator $A\in\b$ is {\em polaroid} if the isolated points of the spectrum of $A$, 
points $\in\iso\sigma(A)$, are poles of (the resolvent of) $A$. It is known, \cite[Theorem 2.4]{D}, that contractive (more generally, power bounded) $m$-isometric Banach 
space operators $S$ (i.e., contractions, respectively power bounded,  $S\in\b$ such that  $\sum_{j=0}^m{(-1)^{j}\left(\begin{array}{clcr}m\\j\end{array}\right)||S^{m-j}x||^2}=0$ for all $x\in\X$) are isometric, hence polaroid. This result extends to power bounded $S, T\in\b$ such that $\triangle^m_{T,S}(I)=0$. We prove in the following that the $n$-th power (hence the operator itself) of an  $n$-quasi $m$-isometric operator in $\B$ is polaroid whenever it is a contraction (more generally, power bounded). Indeed, we prove more: Power bounded operators $S, T\in\B$ such that $[S, T^*]=0$ and $\triangle^m_{T,S}(S^{*n}S^n)=0$ are polaroid.

\

The plan of this paper is as follows. We introduce our notation/terminology, alongwith some complementary results, in Section 2. Here we have a first look at the structure of $n$-quasi $[m,d]$-operators. Section 3 is devoted to proving the polaroid property for $n$-quasi left $m$-invertible operators, Section 4 considers the product of an $n$-quasi $[m_1,d]$-operator with an $[m_2,d]$-operator and Section 5 deals with 
perturbation by nilpotents. As we point out at various points in the paper, our results represent a considerable improvement upon various extant results.

\section {\sfstp Complementary results} Given a Banach space operator $A\in\b$, we denote the isolated points of the spectrum $\sigma(A)$ (resp., the approximate point spectrum $\sigma_a(A)$, the surjectivity spectrum $\sigma_{su}(A)$) of $A$ by $\iso\sigma(A)$ (resp., $\iso\sigma_a(A)$, $\i\sigma_{su}(A)$). Let $A-\lambda$ denote $A-\lambda I$. The operator $A$ is said to have SVEP, {\em the single--valued extension property}, at a point $\lambda$ of the complex plane $\C$ if, for every neighborhood ${\cal O}_{\lambda}$ of $\lambda$, the only analytic function $f:{\cal O}_{\lambda}\longrightarrow\X$ satisfying $(A-\mu)f(\mu)=0$ for all $\mu\in{\cal O}_{\lambda}$
is the function $f\equiv 0$; we say that {\em $A$ has SVEP} if it has SVEP at every $\lambda\in\C$. The  ascent $\asc(A)$ (resp., descent $\dsc(A)$) of $A$ is the least non--negative integer $n$ such that $A^{-n}(0)=A^{-(n+1)}(0)$ (resp., $A^n{\cal X}={A^{n+1}}{\cal X}$); if no such integer exists, then $\asc(A)=\infty$ (resp., $\dsc(A)=\infty$). It is well known, \cite{{A},{He},{LN},{TL}}, that $\asc(A)<\infty$ implies $A$ has SVEP at $0$ and $\dsc(A)<\infty$ implies $A^*$, the dual operator,  has SVEP at $0$, and that finite ascent and descent imply their equality. A point $\lambda\in\iso\sigma(A)$ is a pole of (the resolvent of) $A$ if
$\asc(A-\lambda)=\dsc(A-\lambda)<\infty$.

\

For a given operator $A\in\b$, let {\em $\Pi_a(A)=\{\lambda\in\iso\sigma_a(A):$there exits an integer $d\geq 1$ such that $\asc(A-\lambda)\leq d$ and $(A-\lambda)^{d+1}$ is closed$\}=$ set of left poles of $A$}, and let {\em $\Pi(A)=\{\lambda\in\iso\sigma(A):\asc(A-\lambda)=\dsc(A-\lambda)<\infty\}=$ set of poles of $A$}. Then $\Pi(A)\subseteq\Pi_a(A)$ and a necessary and sufficient condition for $\lambda\in\Pi_a(A)$ to imply $\lambda\in\Pi(A)$ is that $A^*$ has SVEP at $\lambda$ \cite{A}. {\em
	We say that  $A$ is polaroid (resp. left polaroid)} if $\{\lambda\in\sigma(A):\lambda\in\iso\sigma(A)\}=\Pi(A)$ (resp.,
$\{\lambda\in\sigma(A):\lambda\in\iso\sigma_a(A)\}=\Pi_a(A)$.) To every $\lambda\in\iso\sigma(A)$, there corresponds a decomposition $$\X= H_0(A-\lambda)\oplus K(A-\lambda),$$ where $H_0(A-\lambda)$, the quasinilpotent part of $A-\lambda$, and $K(A-\lambda)$, the  analytic core of $A-\lambda$, are the sets $$H_0(A-\lambda)=\{x\in\X:\lim_{n\longrightarrow\infty}{||(A-\lambda)^nx||^{\frac{1}{n}}=0}\}$$ and
\begin{center}$K(A-\lambda)=\{x\in\X:\mbox{there exists a sequence}\hspace{2mm}\{x_n\}\subset\X\hspace{2mm}\mbox{and}\hspace{2mm} \delta>0\hspace{2mm}\mbox{for which}$\\ $x=x_0, (A-\lambda)x_{n+1}=x_n\hspace{2mm}\mbox{and}\hspace{2mm}\hspace{2mm}||x_n||\leq \delta^n||x||\hspace{2mm}\mbox{for all}\hspace{2mm}n=1,2, ...\}$\end{center}
\cite{A}. $H_0(A-\lambda)$ and $K(A-\lambda)$ are generally non-closed hyperinvariant subspaces of $A-\lambda$ such that $(A-\lambda)^{-p}(0)\subseteq H_0(A-\lambda)$ for all positive integers $p$ and $(A-\lambda)K(A-\lambda)=K(A-\lambda)$. A necessary and sufficient condition for  a $\lambda\in\iso\sigma(A)$
to be a pole of $A$ is that $H_0(A-\lambda)=(A-\lambda)^{-n}(0)$ for some integer $n>0$. (The number $n$ is then said to be the order of the pole  at $\lambda$; if $n=1$, then the pole is said to be a simple pole.)

\

Similarities preserve spectrum (hence, isolated points of the spectrum), the ascent and the descent. Hence:
{\em Similarities preserve the polaroid property}. Recall that an $A\in\b$ is an isometry if $||Ax||=||x||$
for all $x\in\X$. Isometries are normaloid operators, i.e., if an $A\in\b$ is isometric then $||A||$ equals the spectral radius
$r(A)=\lim_{n\longrightarrow\infty}{||A^n||^{\frac{1}{n}}}$. The inverse of an isometry, whenever it exists as a bounded operator,
is again an isometry. Since the restriction of an isometry to an invariant subspace is again an isometry, isometries are {\em totally hereditarily normaloid operators} (see \cite{D0}). Conclusion: {\em Invertible isometries are polaroid} (\cite{D0}; see also\cite[Theorem 1.5.13]{LN}).

\

Given operators $S,T\in\b$, it is seen that $$\triangle^{m+k}_{T,S}(I)=(L_TR_S-I)^k(\triangle^m_{T,S}(I))=\sum_{j=0}^k{(-1)^{j}\left(\begin{array}{clcr}k\\j\end{array}\right)T^{k-j}\triangle^m_{T,S}(I)S^{k-j}}$$  and  $$\delta^{m+k}_{T,S}(I)=(L_T-R_S)^k(\delta^m_{T,S}(I))=\sum_{j=0}^k{(-1)^{j}\left(\begin{array}{clcr}k\\j\end{array}\right)T^{k-j}\delta^m_{T,S}(I)S^j}$$ for all integers $m,k\geq 1$. Hence: \begin{lem}\label{lem01} If $d^m_{T,S}(I)=0$, then $d^t_{S,T}(I)=0$ for all integers $t\geq m$.\end{lem}

 For an operator $S\in\B$, let $\overline{S^n(\H)}$  denote the closure of the range of $S^n$, and let ${S^*}^{-n}(0)$ denote the kernel
of ${S^*}^n$. If an operator $T\in\B$ is such that  $[S,T^*]=ST^*-T^*S=0$, then $\H$ has a direct sum decomposition $\H=\overline{S^n(\H)}\oplus {S^*}^{-n}(0)$, and $S, T^*$ have upper triangular representations
\begin{eqnarray}\label{arr0}  S=\left(\begin{array}{clcr}S_1 & S_0\\0 & S_2\end{array}\right), \hspace{4mm} T^*=\left(\begin{array}{clcr}T^*_1 & T^*_0\\0 & T^*_2\end{array}\right),\end{eqnarray}
where $$S^n_2=0 \hspace{3mm}\mbox{and}\hspace{3mm} [S_1,T^*_1]=0.$$
The hypothesis $S^{*n}d^m_{T,S}(I)S^n=0$  implies that if $d=\triangle$ then
\begin{eqnarray*}
	& &  S^{*n}\triangle^m_{T,S}(I)S^n=0 \Longleftrightarrow S^{*n}\{\sum_{j=0}^m{(-1)^{j}\left(\begin{array}{clcr}m\\j\end{array}\right)T^{m-j}S^{m-j}}\}S^n=0\\ &\Longleftrightarrow &
	\left(\begin{array}{clcr} S^{*n}_1 & 0\\X^* & 0\end{array}\right)\{\sum_{j=0}^m{(-1)^{j}\left(\begin{array}{clcr}m\\j\end{array}\right)\left(\begin{array}{clcr} T_1^{m-j}S_1^{m-j} & X_{1j}\\X_{2j} & X_{3j}\end{array}\right)}\}\left(\begin{array}{clcr} S_1^n & X\\0 & 0\end{array}\right)=0
\end{eqnarray*}
and if $d=\delta$ then
\begin{eqnarray*}
	& &  S^{*n}\delta^m_{T,S}(I)S^n=0 \Longleftrightarrow S^{*n}\{\sum_{j=0}^m{(-1)^{j}\left(\begin{array}{clcr}m\\j\end{array}\right)T^{m-j}S^j}\}S^n=0\\ &\Longleftrightarrow &
	\left(\begin{array}{clcr} S^{*n}_1 & 0\\X^* & 0\end{array}\right)\{\sum_{j=0}^m{(-1)^{j}\left(\begin{array}{clcr}m\\j\end{array}\right)\left(\begin{array}{clcr} T_1^{m-j}S_1^j & X_{1j}\\X_{2j} & X_{3j}\end{array}\right)}\}\left(\begin{array}{clcr} S_1^n & X\\0 & 0\end{array}\right)=0
\end{eqnarray*}
for some operators $X$ and $X_{ij}$ ($i=1,2,3$). Hence
\begin{eqnarray*} & &
	S^{*n}_1\{\sum_{j=0}^m{(-1)^{j}\left(\begin{array}{clcr}m\\j\end{array}\right)T_1^{m-j}S_1^{m-j}}\}S^n_1=0\\
	&\Longleftrightarrow& \sum_{j=0}^m{(-1)^{j}\left(\begin{array}{clcr}m\\j\end{array}\right)T_1^{m-j}S_1^{m-j}}=0,\end{eqnarray*} and
\begin{eqnarray*} & &
	S^{*n}_1\{\sum_{j=0}^m{(-1)^{j}\left(\begin{array}{clcr}m\\j\end{array}\right)T_1^{m-j}S_1^j}\}S^n_1=0\\
	&\Longleftrightarrow& \sum_{j=0}^m{(-1)^{j}\left(\begin{array}{clcr}m\\j\end{array}\right)T_1^{m-j}S_1^j}=0,\end{eqnarray*}
i.e., $d^m_{T_1,S_1}(I)=0$.  Consequently, \cite[Remark 2.7]{DM} and Lemma \ref{lem01}, $d^m_{T_1^p,S_1^p}(I)=0$ for every integer $p\geq 1$. Hence $$S^{*n}d^m_{T^p,S^p}(I)S^n=0, \hspace{2mm}\mbox{for all integers}\hspace{2mm} p\geq 1.$$ The observations that
$$\triangle^{m+1}_{A,B}(I)=A\triangle^m_{A,B}(I) B-\triangle^m_{A,B}(I), \hspace{2mm} \delta^{m+1}_{A,B}(I)=A\delta^m_{A,B}(I)-\delta^m_{A,B}(I) B $$ lead to the implication $$d^m_{T_1,S_1}(I_1)=0\Longleftrightarrow d^t_{T_1,S_1}(I_1)=0 \hspace{2mm}\mbox{for all integers}\hspace{2mm} t\geq m,$$
and hence $$S^{*n}d^m_{T,S}(I)S^n=0\Longrightarrow S^{*n}d^t_{T,S}(I)S^n=0 \hspace{2mm}\mbox{for all integers}\hspace{2mm} t\geq m.$$  If we let $X$ denote the operator
$$X=\sum_{j=0}^{n-1}{S_1^{n-1-j}S_0S^j_2},$$then
$$S^n=\left(\begin{array}{clcr}S^n_1 & X\\0 & 0\end{array}\right).$$ Now if $S^{*n}\triangle^m_{T,S}(I)S^n=0$, then $\triangle^m_{T_1,S_1}(I_1)=0$ implies $S_1$ is ($m$-left invertible, hence) left invertible. Consequently, if $S_1$ has a dense range (or, equivalently, $S_1^*$ has SVEP at $0$), then the operator $S^n$ is similar to  $A=S^n_1\oplus 0$   (with the similarity implemented by the invertible operator $E=\left(\begin{array}{clcr}S^n_1 & S^n_1X\\0 & 1\end{array}\right)$). Observe that the operator $A$ is not left $m$-invertible (i.e., there does not exist an operator $B\in\B$ such that $\triangle^m_{B,A}(I)=0$). Letting  $T=S^*$ (so that $S^{*n}\triangle^m_{S^*,S}(I)S^n=0$ -- such operators have been called {\em $n$-quasi $m$-isometric} [11]), it then follows that  $S^n_1$ is $m$-isometric and, if $S_1$ has a dense range, $S^n$ is similar to $A$. Operators $S\in\B$ for which $\delta^m_{S^*,S}(I)=0$ are called {\em $m$-selfadjoint operators} [7]. If $S^{*n}\delta^m_{S^*,S}(I)S^n=0$ (i.e., if $S$ is {\em $n$-quasi $m$-selfadjoint}), then  ($S_1$, hence) $S^p_1$ is $m$-selfadjoint for all integers $p\geq 1$ [7]. More is true, as we prove in the following.

\vskip4pt

Given a positive operator $(0\leq) Q\in\B$, we say that the operator $S\in\B$ is $[m, Q]$-isometric (resp., $[m, Q]$-selfadjoint) if $\triangle^m_{S^*,S}(Q)=0$ (resp., $\delta^m_{S^*,S}(Q)=0$); we say that $S\in [m, d(Q)]$  if $d^m_{S^*,S}(Q)=0$, $d=\triangle$ or $\delta$. We assume in the following that $S^n_1=(S|_{\overline{S^n(\H)}})^n$ has the polar decomposition $S^n_1=U_1P_1$. It is then clear that $U_1$ is an isometry and $P_1\geq 0$ is invertible in the case in which $S$ is $n$-quasi $m$-isometric, and $U_1$ is isometric and $P_1\geq 0$ is injective in the case in which $S$ is $n$-quasi $m$-selfadjoint and injective.  Define the operator $P\in B({\overline{S^n(\H)}}\oplus S^{*{-n}}(0))$ by $P=P_1\oplus I_2$.

\begin{pro}\label{pro00}   Let $S\in\B$ be such that $S^{*n}d^m_{S^*,S}(I)S^n=0$ for some integers $m,n\geq 1$.

\vskip4pt\noindent (i) If $d=\triangle$, then there exist operators $Q,A\in\B$ such that $Q\geq 0$, $\triangle^m_{A^*,A}(Q)=0$ and $S^n$ is similar to $A$.

\vskip4pt\noindent (ii) If $d=\delta$ and the operator $S$ is injective, then  there exist operators $Q,A\in\B$ such that $Q\geq 0$, $\delta^m_{A^*,A}(Q)=0$ and $\delta_{A,S^n}(P)=0$.

\vskip4pt\noindent (iii) If $S$ is left  invertible, then there exists an operator $B\in\B$ such that $d^m_{B^*,B}(I)=0$ and $S^n$ is similar to $B$.
\end {pro}
\begin{demo}    The hypothesis $S^{*n}d^m_{S^*,S}(I)S^n=0$   implies $d^m_{S^{*p}_1,S_1^p}(I)=0$, and hence $$S^{*n}d^m_{S^{*p},S^p}(I)S^n=0 \hspace{2mm}\mbox{for all integers}\hspace{2mm} p\geq 1.$$
	Let, as above, $$S^n=\left(\begin{array}{clcr}S^n_1 & X\\0 & 0\end{array}\right).$$  Define the operators $A_1$, $A$ and $Q$ by
	$$A_1=P_1U_1,\hspace{2mm}  A=\left(\begin{array}{clcr}A_1 & P_1X\\0 & 0\end{array}\right)  \hspace{2mm} \mbox{and}\hspace{2mm} \hspace{2mm} \mbox{and}\hspace{2mm} Q=\left(\begin{array}{clcr}I_1 & U^*_1X\\X^*U_1 & X^*X\end{array}\right).$$ Let $I_2$ denote (as above)  the identity of $B(S^{-*n}(0))$.

\vskip4pt\noindent (i). If $d=\triangle$, then (upon letting $p=n$ in the above) we have:

	\begin{eqnarray*} & &  S^{*n}\triangle^m_{S^{*n},S^n}(I)S^n=0 \Longleftrightarrow
	 S^{*n}\{\sum_{j=0}^m{(-1)^{j}\left(\begin{array}{clcr}m\\j\end{array}\right){S^*}^{n(m-j)}S^{n(m-j)}}\}S^{n}=0\\ &\Longleftrightarrow& (P_1\oplus I_2)\left(\begin{array}{clcr}U_1^* & 0\\X^* & 0\end{array}\right)\{\sum_{j=0}^m{(-1)^{j}\left(\begin{array}{clcr}m\\j\end{array}\right)\left(\begin{array}{clcr}S^{*n}_1 & 0\\X^* & 0\end{array}\right)^{m-j} \left(\begin{array}{clcr}S^n_1 & X\\0 & 0\end{array}\right)^{m-j}}\}\times\\ & \times&\left(\begin{array}{clcr}U_1 & X\\0 & 0\end{array}\right)(P_1\oplus I_2)=0\\  &\Longleftrightarrow& \sum_{j=0}^m{(-1)^{j}\left(\begin{array}{clcr}m\\j\end{array}\right)\left(\begin{array}{clcr}U_1^* & 0\\X^* & 0\end{array}\right) \left(\begin{array}{clcr}P_1U_1^* & 0\\X^* & 0\end{array}\right)^{m-j}\left(\begin{array}{clcr}U_1P_1 & X\\0 & 0\end{array}\right)^{m-j}\left(\begin{array}{clcr}U_1 & X\\0 & 0\end{array}\right)}=0\\ &\Longleftrightarrow& \sum_{j=0}^m{(-1)^{j}\left(\begin{array}{clcr}m\\j\end{array}\right)\left(\begin{array}{clcr}A_1^* & 0\\X^*P_1 & 0\end{array}\right)^{m-j}}{\left(\begin{array}{clcr}U_1^* & 0\\X^* & 0\end{array}\right)\left(\begin{array}{clcr}U_1 & X\\0 & 0\end{array}\right)\left(\begin{array}{clcr}A_1 & P_1X\\0 & 0\end{array}\right)^{m-j}}=0\\ &\Longleftrightarrow& \sum_{j=0}^m{(-1)^{j}\left(\begin{array}{clcr}m\\j\end{array}\right)\left(\begin{array}{clcr}A^*_1 & 0\\X^*P_1 & 0\end{array}\right)^{m-j}\left(\begin{array}{clcr}I_1 & U_1^*X\\X^*U_1 & X^*X\end{array}\right)\left(\begin{array}{clcr}A_1 & P_1X\\0 & 0\end{array}\right)^{m-j}}=0\\ &\Longleftrightarrow& \sum_{j=0}^m{(-1)^{j}\left(\begin{array}{clcr}m\\j\end{array}\right)\left(\begin{array}{clcr}A_1^* & 0\\X^*P_1 & 0\end{array}\right)^{m-j} Q \left(\begin{array}{clcr}A_1 & P_1X\\0 & 0\end{array}\right)^{m-j}}=0\\ &\Longleftrightarrow&
		\sum_{j=0}^m{(-1)^{j}\left(\begin{array}{clcr}m\\j\end{array}\right)A^{*(m-j)} Q A^{m-j}}=0\\ &\Longleftrightarrow& \triangle^m_{A^*,A}(Q)=0.  \end{eqnarray*}
	 Set $P_1\oplus I_2=P$. Then  $$S^n=\left(\begin{array}{clcr}U_1P_1 & X\\0 & 0\end{array}\right)=P^{-1}\left(\begin{array}{clcr}A_1 & P_1X\\0 & 0\end{array}\right)P=P^{-1}AP,$$ i.e., $S^n$ is similar to $A$.

\vskip4pt\noindent (ii). If $d=\delta$, then (following the notation developed above): \begin{eqnarray*} & &  S^{*n}\delta^m_{S^*,S}(I)S^n=0 \Longrightarrow   S^{*n}\delta^m_{S^{*n},S^n}(I)S^n=0\\  &\Longleftrightarrow& \sum_{j=0}^m{(-1)^{j}\left(\begin{array}{clcr}m\\j\end{array}\right) S^{*n(m-j+1)}S^{n(j+1)}}=0\\ &\Longleftrightarrow& \sum_{j=0}^m{(-1)^{j}\left(\begin{array}{clcr}m\\j\end{array}\right) S^{*n(m-j)}(P_1\oplus I_2)\left(\begin{array}{clcr}U_1^* & 0\\X^* & 0\end{array}\right)   \left(\begin{array}{clcr}U_1 & X\\0 & 0\end{array}\right)(P_1\oplus I_2)S^{nj}}=0\\ &\Longleftrightarrow&
P\{\sum_{j=0}^m{(-1)^{j}\left(\begin{array}{clcr}m\\j\end{array}\right) \left(\begin{array}{clcr}U_1^*P_1 & 0\\X^*P_1 & 0\end{array}\right) ^{m-j}\left(\begin{array}{clcr}I_1 & U_1^*X\\X^*U_1 & X^*X\end{array}\right) \left(\begin{array}{clcr}P_1U_1 & P_1X\\0 & 0\end{array}\right)^j}\}P=0\\ &\Longleftrightarrow&
		\sum_{j=0}^m{(-1)^{j}\left(\begin{array}{clcr}m\\j\end{array}\right)A^{*(m-j)} Q A^j}=0\\ &\Longleftrightarrow& \delta^m_{A^*,A}(Q)=0.\end{eqnarray*}
Evidently, $P=P_1\oplus I_2$ is a quasi-affinity such that $\delta_{A,S^n}(P)=0$.

\vskip4pt\noindent (iii). Assume now that $S$ is left invertible. Then $P$ and $Q$ (defined as above) are positive invertible, and \begin{eqnarray*} & & S^{*n}\triangle^m_{S^*,S}(I)S^n=0\Longrightarrow S^{*n}\triangle^m_{S^{*n},S^n}(I)S^n=0\\ &\Longrightarrow& \sum_{j=0}^m{(-1)^{j}\left(\begin{array}{clcr}m\\j\end{array}\right)A^{*(m-j)}QA^{m-j}}=0\\ &\Longleftrightarrow& \sum_{j=0}^m{(-1)^{j}\left(\begin{array}{clcr}m\\j\end{array}\right)(Q^{-{\frac{1}{2}}}A^*Q^{{\frac{1}{2}}})^{m-j} (Q^{{\frac{1}{2}}}AQ^{-{\frac{1}{2}}})^{m-j}}=0\end{eqnarray*} and
\begin{eqnarray*} & & S^{*n}\delta^m_{S^*,S}(I)S^n=0\Longrightarrow S^{*n}\delta^m_{S^{*n}S^n}(I)S^n=0\\ &\Longrightarrow& \sum_{j=0}^m{(-1)^{j}\left(\begin{array}{clcr}m\\j\end{array}\right)A^{*(m-j)}QA^{j}}=0\\ &\Longleftrightarrow& \sum_{j=0}^m{(-1)^{j}\left(\begin{array}{clcr}m\\j\end{array}\right) (Q^{-{\frac{1}{2}}}A^*Q^{{\frac{1}{2}}})^{m-j} (Q^{{\frac{1}{2}}}AQ^{-{\frac{1}{2}}})^j}=0.\end{eqnarray*}  Now define $B\in\B$ by $$B=Q^{{\frac{1}{2}}}AQ^{-{\frac{1}{2}}};$$ then $$d^m_{B^*,B}(I)=0.$$  Since $$B=Q^{{\frac{1}{2}}}AQ^{-{\frac{1}{2}}}=Q^{{\frac{1}{2}}}PS^nP^{-1}Q^{-{\frac{1}{2}}}=LS^nL^{-1}\Longrightarrow S^n=L^{-1}BL, \hspace{2mm} L=Q^{\frac{1}{2}}P,$$ $S^n$ is similar to $B$.
\end{demo}
Let $\D$ denote the open unit disc in $\C$ and let $\partial{\D}$ denote the boundary of $\D$.
\begin{cor}\label{cor00}  ({\it cf.}\cite[Corollary 4.3]{SS}) If $d=\triangle$ in the statement of  Proposition \ref{pro00} and the operator $Q$ (in the proof of the proposition) is injective, then $\sigma_p(S)\subseteq \partial{\D}$.\end{cor} \begin{demo} The hypotheses imply $\sum_{j=0}^m{(-1)^{j}\left(\begin{array}{clcr}m\\j\end{array}\right) ||Q^{\frac{1}{2}}A^{m-j}x||^2}=0$ for all $x\in\H$. Consider a $\lambda\in\sigma_p(S)$ such that $Ax=\lambda x$. Then, since $Q$ is injective,  \begin{eqnarray*} & & \sum_{j=0}^m{(-1)^{j}\left(\begin{array}{clcr}m\\j\end{array}\right) |\lambda|^{2(m-j)} ||Q^{\frac{1}{2}}x||^2}=0\\ &\Longleftrightarrow& (1-|\lambda|^2)^m=0 \Longleftrightarrow |\lambda|=1.\end{eqnarray*} Since $S^n$ is similar to $A$, $\sigma_p(S)^n=\sigma_p(S^n)=\sigma_p(A)\subseteq \partial(\D)$.\end{demo}

 Proposition \ref{pro00} is a generalization of some extant results. For example, if $d=\triangle$, $n=1$ and $m=2$, then $S^*\triangle^2_{S^*,S}(I)S=0$ (i.e., $S$ is $1$-quasi $2$-isometric) implies $\triangle^2_{A^*,A}(Q)=0$ (where the operators $A, Q$ are as defined in the proof of the proposition and the operator $S$ is similar to $A$); if also $S$ is left invertible, then $\triangle^2_{B^*,B}(I)=0$ (i.e., $B$ is $2$-isometric) for some operator $B$ similar to the operator $S$  ({\it{cf.}} \cite[Theorem 2.5]{MP}).   In their considerations on the spectral properties of  $A$-contractions, L. Suciu and N. Suciu \cite{SS} define an operator $S\in\B$ to be $n$-quasi isometric if $S^{*n}(S^*S-1)S^n=0$. In our terminology this equates to  $S^{*n}\triangle_{S^*,S}(I)S^n=0$ (equivalently, ``$S$ is $n$-quasi $1$-isometric"). Thus, for  $n$-quasi isometric operators $S$, $S^n_1$ is isometric; indeed, since $S^{*n}_1(S^*_1S_1-I_1)S^n_1=0$, $S_1$ is isometric. Assume now that $n=1$ and $0$ is a normal eigenvalues of $S$ (i.e., $S^{-1}(0)\subseteq S^{*-1}(0)$). Then $S=S_1\oplus 0$ is a partial isometry ({\it cf.} \cite[Theorem 3.12 and Corollary 3.13]{SS}). For a general $n$-quasi isometry $S$, $S=\left(\begin{array}{clcr}S_1&S_0\\0&S_2\end{array}\right)\in B({\overline{S^n(\H)}}\oplus S^{*-n}(0))$, where $S_1$ is isometric and $S_2$ is $n$-nilpotent. Consequently, $S$ has SVEP and hence  \cite[Theorem 4.6]{SS}:  (i) $\sigma(S)=\overline{\sigma_a(S^*)}$.  (ii) $\sigma(S)=\overline{\D}$, the closed unit disc, if $S_1$ is not invertible and $\sigma(S)\subseteq \partial{\D}\cup\{0\}$ if $S_1$ is invertible. In either case, $\sigma_a(S)\subseteq \partial{\D}\cup\{0\}$. (iii) If $\lambda, \mu$ are two distinct non-zero eigenvales of $S$,  then $\lambda, \mu\in\sigma_p(S_1)$ and the corresponding eigenspaces are mutually orthogonal. Observe that if $n=1$, then $S_1$ is isometric. If also $||S||\leq 1$, then $S^pS^{*p}=S^p_1S^{*p}_1+S^{p-1}_1S_0S^*_0S^{*(p-1)}_1\oplus 0$ is a contraction (thus: $S^p_1S^{*p}_1+S^{p-1}_1S_0S^*_0S^{*(p-1)}_1\leq I_1$). Consequently, $$S^{*p}S^p=\left(\begin{array}{clcr} I_1&S^{*p}_1S_0\\S_0^*S^{p-1}_1&S_0^*S_0\end{array}\right)\geq S^pS^{*p}$$ for all integers $p\geq 1$ \cite[Theorem 3.3]{SS}.

\vskip6pt \noindent

 Let $C$ be a conjugation of $\H$ (i.e., $C:\H\longrightarrow\H$ is a conjugate-linear operator such that $C^2=I$ and $<Cx,y>=<Cy,x>$ for all $x,y\in\H$).  If one chooses $T=CS^*C$ in $d^m_{T,S}(I)=0$, then $$\triangle^m_{CS^*C,S}(I)=0\Longleftrightarrow \triangle^m_{S^*,CSC}(I)=0$$ defines the class of {\em $(m,C)$-isometric operators} and $$\delta^m_{CS^*C,S}(I)=0\Longleftrightarrow \delta^m_{S^*,CSC}(I)=0$$ defines the class of {\em $(m,C)$-symmetric operators} \cite{{CLM}, {CKL}}. It is known \cite{{CLM}, {CKL}} that $$d^m_{S^*,CSC}(I)=0\Longleftrightarrow d^t_{S^*,CSC}(I)=0 \hspace{2mm}\mbox{for all integers}\hspace{2mm} t\geq m$$
and $$d^m_{S^*,CSC}(I)=0\Longleftrightarrow d^m_{S^{*p},CS^pC}(I)=0 \hspace{2mm}\mbox{for all integers}\hspace{2mm} p\geq 1.$$

\vskip3pt

It is clear that if $S^{*n}d^m_{S^*,CSC}(I)S^n=0$, then $S\in B({\overline{S^n(\H)}}\oplus S^{*{-n}}(\H))$ has a representation
$$CS^nC= C\left(\begin{array}{clcr}S^n_1 & X\\0 & 0\end{array}\right)C= C\left(\begin{array}{clcr} U_1P_1 & X\\0 & 0\end{array}\right)C $$
(where the operator $X$ is as defined above). In particular, if the conjugation $C:  {\overline{S^n(\H)}}\oplus S^{*{-n}}(\H)\longrightarrow {\overline{S^n(\H)}}\oplus S^{*{-n}}(\H)$ has a representation $C=C_1\oplus C_2$, then $$d^m_{S^*_1,C_1S_1C_1}(I_1)=0\Longrightarrow S^{*n}d^m_{S^{*p},CS^pC}(I)S^n=0$$ for all integers $p\geq 1$. If, now, $S$ satisfies the additional property that $CSCS=S^2$, then $$S^{*n}d^m_{S^*,CSC}(I)S^n=0\Longleftrightarrow S^{*n}d^m_{S^*,S}(I)S^n=0$$ and Proposition \ref{pro00} applies. In general,  Proposition \ref{pro00} seemingly does not extend to operators $S$ satisfying $S^{*n}d^m_{S^*,CSC}(I)S^n=0$. Define the operator $M\in\B$ by $$M=\left(\begin{array}{clcr} U_1 & X\\ 0 & 0\end{array}\right)$$ (where $U_1$ and $X$ are the operators defined above). The following proposition says that a result very similar to Proposition \ref{pro00} holds in the case in which $[C,M]=0$ and  $C=C_1\oplus C_2$.

\begin{pro}
		\label{pro110} Let $S\in\B$ be such that $S^{*n}d^m_{S^*,CSC}(I)S^n=0$ (so that $S$ is either $n$-quasi $(m,C)$-isometric or $S$ is $n$-quasi $(m,C)$ symmetric), where the conjugation  $C=C_1\oplus C_2: \overline{S^n(\H)}\oplus {S^*}^{-n}(0)\longrightarrow \overline{S^n(\H)}\oplus {S^*}^{-n}(0)$  satisfies $[C,M]=0$.
\vskip4pt\noindent (i) If $d=\triangle$, then there exist operators $Q, A\in\B$ such that $Q\geq 0$, $\triangle^m_{A^*,CAC}(Q)=0=\triangle^m_{CA^*C,A}(CQC)$ and $S^n$ is similar to $A$.

\vskip4pt\noindent (ii) If $d=\delta$ and the operator $S$ is injective, then  there exist operators $Q,A\in\B$ such that $Q\geq 0$, $\delta^m_{A^*,CAC}(Q)=0=\delta^m_{CA^*C,A}(CQC)$ and $\delta_{A,S^n}(P)=0$.

\vskip4pt\noindent (iii) If $S$ is left invertible, then there exists an operator $B\in\B$ such that $d^m_{B^*,CBC}(I)=0$ and $S^n$ is similar to $B$.    \end{pro}

\begin{demo}  We start by observing that  $$S^{*n}d^m_{S^*,CSC}(I)S^n=0\Longleftrightarrow S^{*n}d^m_{S^{*p},CS^pC}(I)S^n=0$$ for all integers $p\geq 1$.

\vskip4pt\noindent (i).  Case $d=\triangle$.  Following the notation of the proof of Proposition \ref{pro00}, we have:
\begin{eqnarray*} & & S^{*n}d^m_{S^{*n},CS^nC}(I)S^n=0\\ &\Longleftrightarrow& S^{*n}\{\sum_{j=0}^m{(-1)^{j}\left(\begin{array}{clcr}m\\j\end{array}\right)S^{*n(m-j)}CS^{n(m-j)}C}\}S^n=0\\ &\Longleftrightarrow& \sum_{j=0}^m{(-1)^{j}\left(\begin{array}{clcr}m\\j\end{array}\right){S^*}^{n(m-j+1)}CS^{n(m-j+1)}C}=0\\ &\Longleftrightarrow& \sum_{j=0}^m {(-1)^{j}\left(\begin{array}{clcr}m\\j\end{array}\right)\left(\begin{array}{clcr}P_1U^*_1 & 0\\X^* & 0\end{array}\right)^{m-j+1} C \left(\begin{array}{clcr}U_1P_1 & X\\0 & 0\end{array}\right)^{m-j+1}C}=0 \\  &\Longleftrightarrow& \left(\begin{array}{clcr}P_1 & 0\\0 & I_2\end{array}\right) \{ {\sum_{j=0}^m{(-1)^{j}\left(\begin{array}{clcr}m\\j\end{array}\right)\left(\begin{array}{clcr}U^*_1P_1 &
			0\\X^*P_1 & 0\end{array}\right)^{m-j} \left(\begin{array}{clcr}U^*_1 & 0\\X^*& 0\end{array}\right)\times}}\\ &\times& {\left(\begin{array}{clcr}C_1U_1C_1 &
		C_1XC_2\\0 & 0\end{array}\right)\left(\begin{array}{clcr}C_1P_1U_1C_1 & C_1P_1XC_2\\0 & 0\end{array}\right)^{m-j} }\} {\left(\begin{array}{clcr}C_1P_1C_1 & 0\\0 & I_2\end{array}\right)=0} \\ &\Longleftrightarrow& \sum_{j=0}^m{(-1)^{j}\left(\begin{array}{clcr}m\\j\end{array}\right)\left(\begin{array}{clcr}U^*_1P_1 &
			0\\X^*P_1 & 0\end{array}\right)^{m-j} \left(\begin{array}{clcr}I_1 & U^*_1C_1XC_2\\X^*U_1C_1U_1 & X^*C_1XC_2\end{array}\right)\times}\\ &\times& {\left(\begin{array}{clcr}C_1P_1U_1C_1 & C_1P_1XC_2\\0 & 0\end{array}\right)^{m-j} }=0 .\end{eqnarray*}

By hypothesis, $[C,M]=0$. Hence $$C_1X=XC_2, [C_1,U_1]=0\hspace{2mm}\mbox{and}\hspace{2mm} \left(\begin{array}{clcr}I_1 & U^*_1C_1XC_2\\X^*U_1C_1U_1 & X^*C_1XC_2\end{array}\right)=\left(\begin{array}{clcr}I_1 & U^*_1X\\X^*U_1 & X^*X\end{array}\right)=Q$$ for some positive operator $Q$. Consequently \begin{eqnarray*} & &  S^{*n}\triangle^m_{S^*,CSC}(I)S^n=0\\ &\Longrightarrow&	 \sum_{j=0}^m{(-1)^{j}\left(\begin{array}{clcr}m\\j\end{array}\right)\left(\begin{array}{clcr}U^*_1P_1 & 0\\X^*P_1 & 0\end{array}\right)^{m-j}  Q \left(\begin{array}{clcr}C_1P_1U_1C_1 & C_1P_1XC_2\\0 & 0\end{array}\right)^{m-j} }=0\\ &\Longrightarrow& \triangle^m_{A^*,CAC}(Q)=0\Longleftrightarrow \triangle^m_{CA^*C,A}(CQC)=0,\end{eqnarray*} where, as before, the operator $A$ is defined by $A=\left(\begin{array}{clcr} P_1U_1&  P_1X\\0 & 0\end{array}\right)= PS^nP^{-1}$.

\vskip4pt\noindent (ii). Case $d=\delta$. The hypothesis $S$ is injective implies $P\geq 0$ has a dense range. Using the same notation as above, we have:
\begin{eqnarray*} & &  S^{*n}\delta^m_{S^*,CSC}(I)S^n=0\\ &\Longrightarrow&  \sum_{j=0}^m{(-1)^{j}\left(\begin{array}{clcr}m\\j\end{array}\right){S^*}^{n(m-j+1)}CS^{n(j+1)}C}=0\\ &\Longleftrightarrow&
\left(\begin{array}{clcr}P_1 & 0\\0 & I_2\end{array}\right)\{\sum_{j=0}^m{(-1)^{j}\left(\begin{array}{clcr}m\\j\end{array}\right)\left(\begin{array}{clcr}U^*_1P_1 &
			0\\X^*P_1 & 0\end{array}\right)^{m-j} \left(\begin{array}{clcr}U^*_1 & 0\\X^*& 0\end{array}\right)}\times \\ & \times& {\left(\begin{array}{clcr}C_1U_1C_1 &
		C_1XC_2\\0 & 0\end{array}\right)\left(\begin{array}{clcr}C_1P_1U_1C_1 & C_1P_1XC_2\\0 & 0\end{array}\right)^j }\}{\left(\begin{array}{clcr}C_1P_1C_1 & 0\\0 & I_2\end{array}\right)=0}\\ &\Longleftrightarrow&
\sum_{j=0}^m{(-1)^{j}\left(\begin{array}{clcr}m\\j\end{array}\right)\left(\begin{array}{clcr}U^*_1P_1 & 0\\X^*P_1 & 0\end{array}\right)^{m-j}  Q \left(\begin{array}{clcr}C_1P_1U_1C_1 & C_1P_1XC_2\\0 & 0\end{array}\right)^j }=0\\ &\Longrightarrow& \sum_{j=0}^m{(-1)^{j}\left(\begin{array}{clcr}m\\j\end{array}\right) A^{*(m-j)}QCA^jC}=0\\ &\Longleftrightarrow&
\delta^m_{A^*,CAC}(Q)=0\Longleftrightarrow \delta^m_{CA^*C,A}(CQC)=0.\end{eqnarray*} It being evident that $\delta_{A,S^n}(P)=0$, $P\geq 0$ a quasi-affinity, the proof is complete.

\vskip4pt\noindent (iii). Arguing as in the proof of Proposition \ref{pro00}, it is seen that $P=P_1\oplus I_2 >0$ and $Q>0$ are invertible; furthermore $CQC=Q$. Since \begin{eqnarray*} & & \sum_{j=0}^m{(-1)^{j}\left(\begin{array}{clcr}m\\j\end{array}\right) A^{*a_j}QCA^{b_j}C}=0\\ &\Longleftrightarrow&
\sum_{j=0}^m{(-1)^{j}\left(\begin{array}{clcr}m\\j\end{array}\right) (Q^{-\frac{1}{2}} A^* Q^{\frac{1}{2}})^{a_j} (CQ^{\frac{1}{2}}A Q^{-\frac{1}{2}}C)^{b_j}}=0\end{eqnarray*} for all positive integers $a_j,$ and $b_j$ we have $$ d^m_{B^*,CBC}(I)=0\Longleftrightarrow d^m_{CB^*C,B}(I)=0;\hspace{2mm} B=Q^{\frac{1}{2}}A Q^{-\frac{1}{2}}$$ Clearly, $S^n=P^{-1}Q^{-\frac{1}{2}}BQ^{\frac{1}{2}}P$ is similar to $B$.
\end{demo}

\section {\sfstp  The Polaroid Property}    If  $\triangle^m_{T,S}(I)=0$ for some $S,T\in\b$ (i.e.,  if $S\in\b$ is left $m$-invertible by $T\in\b)$, then $0\notin\sigma_a(S)$ (for if $0\in\sigma_a(S)$ and $\{x_n\}\subset\X$ is a sequence of unit vectors such that $\lim_{n\longrightarrow\infty}{Sx_n}=0$, then
$$\lim_{n\longrightarrow\infty}||x_n||=\lim_{n\longrightarrow\infty}||\triangle^m_{T,S}(I)x_n||=\lim_{n\longrightarrow\infty}||\sum_{j=0}^m{(-1)^{j}\left(\begin{array}{clcr}m\\j\end{array}\right)T^{m-j}S^{m-j}x_n||}=0$$
-a contradiction). Indeed, if $\lambda\in\sigma_a(S)$, and $\{x_n\}\subset\X$ is a sequence of unit vectors such that $\lim_{n\longrightarrow\infty}(S-\lambda)x_n=0$, then
\begin{eqnarray*} & & \lim_{n\longrightarrow\infty}{\triangle^m_{T,S}(I) x_n}= \lim_{n\longrightarrow\infty}\{\sum_{j=0}^m{(-1)^{j}\left(\begin{array}{clcr}m\\j\end{array}\right)T^{m-j}S^{j} x_n}\}\\ &=&
	\lim_{n\longrightarrow\infty}\{\sum_{j=0}^m{(-1)^{j}\left(\begin{array}{clcr}m\\j\end{array}\right)(\lambda T)^{m-j}}x_n\}\\
	&=& \lim_{n\longrightarrow\infty}(1-\lambda T)^mx_n=0\Longrightarrow \frac{1}{\lambda} \in\sigma_a(T).
\end{eqnarray*}  A similar argument, using this time the fact that $$\triangle^m_{S^*,T^*}(I)=\sum_{j=0}^m{(-1)^{j}\left(\begin{array}{clcr}m\\j\end{array}\right) {S^{*(m-j)}}{T^{*(m-j)}}}=0,$$ shows that $\lambda\in\sigma_{su}(T)$ implies $\frac{1}{\lambda}\in\sigma_{su}(S)$ for all non-zero $\lambda$. (Here $\sigma_{su}(.)$ denotes the surjectivity spectrum.)

\

If we assume $S$, $\triangle^m_{T,S}(I)=0$, to be a contraction  satisfying $\sigma(S)=\overline{\D}$, then $\iso\sigma(S)=\emptyset$ and $S$ is (vacuously) polaroid. If, instead, we assume that $S$ is an invertible contraction with spectrum a subset of the boundary $\partial\D$  of the unit disc $\D$, then $S$ is normaloid (i.e., $||S||=r(S)$) and $\sigma(S)$ consists of the peripheral spectrum ($=\{\lambda:|\lambda|=r(S)\}$) of $S$. The normaloid property of $S$ implies that $\asc(S-\lambda)\leq 1$ and $\dim(\X\setminus{(S-\lambda)(\X)})>0$ \cite[Proposition 54.2]{He}. Thus, if the range $(S-\lambda)^d(\X)$ is closed for some integer $d\geq 1$, then $(S-\lambda)(\X)$ is closed \cite[Proposition 4.10.4]{LN} and $\asc(S-\lambda)\leq 1$, i.e., $\lambda$ is a left pole of $S$. Since $\lambda$ is a boundary point of the spectrum, $\lambda$ is indeed a pole of $S$. Conclusion: ``A necessary and sufficient condition for a point $\lambda\in\sigma(S)$ to be a pole of $S$ for a given left $m$-invertible contraction $S$ (i.e., a contraction $S$ such that $\triangle^m_{T,S}(I)=0$ for some $T\in\b$) with $\sigma(S)\subseteq\partial{\D}$  is that $(S-\lambda)(\X)$ is closed."

\

The hypothesis that $S$ is a left $m$-invertible contraction (resp., $T$ is a right $m$-invertible contraction), even that $S$ is an invertible contraction (resp., $T$ is an invertible  contraction), is not sufficient for $S$ to be polaroid. For example, the operator $S=(I+Q)^{-1}$, $I$ the identity operator and $Q$ the Volterra integration operator, is invertible  with $\sigma(S)=\{1\}$ and $||S||=1$ \cite[Solution 190, Page 302]{Ha}. Since $(I+Q)^{-1}-I=-Q(I+Q)^{-1}=-(I+Q)^{-1}Q$ and $||((I+Q)^{-1}-I)^n)||^{\frac{1}{n}}\leq ||(I+Q)^{-1}||||Q^n||^{\frac{1}{n}}$ converges to $0$ as $n\longrightarrow \infty$, $S$ is not polaroid. Again, if we let $T=(I+Q)^{-1}$ and $S=I+Q$, then $S$ is not polaroid.  A sufficient condition for an operator $S$, $\triangle^m_{T,S}(I)=0$, to be polaroid is that both $S, T$ are power bounded. We recall: $A\in\b$ is power bounded if there exists a positive scalar $M$ such that $\rm{sup_{n\in\n}}||A^n||< M$.
\begin{thm}\label{thm10} If $S, T\in\b$ satisfy $\triangle^m_{T,S}(I)=0$  for some integer $m\geq 1$, then a sufficient condition for $S$ to be polaroid is that $S, T$ are power bounded.	
\end{thm} \begin{demo} If $S, T$ are power bounded, then there exist scalars $M_1, M_2$ such that $$\sup_{n\in\n}||S^n||< M_1, \hspace{2mm} \sup_{n\in\n}||T^n||< M_2$$ (and hence $r(S)=r(T)=1$). This, in view of the fact that ($0\notin\sigma_a(S)$ and) $\{\frac{1}{\lambda}: 0\neq\lambda\in\sigma_a(S)\}\subseteq\sigma_a(T)$ implies $\sigma_a(S)\subseteq\partial{\D}$. Hence
	$$\sigma(S)=\overline{\D}\hspace{2mm}\mbox{if}\hspace{2mm} S\hspace{2mm}\mbox{is not invertible and}\hspace{2mm}\sigma(S)\subseteq\partial{\D} \hspace{2mm}\mbox{if}\hspace{2mm} S\hspace{2mm}\mbox{is invertible}.$$
	Trivially, $S$ is polaroid in the case in which $\sigma(S)=\overline{\D}$. Assume hence that $S$ is invertible (so that $\sigma(S)\subseteq\partial{\D}$). Since $\triangle^m_{T,S}(I)=0$ implies $\triangle^m_{T^p,S^p}(I)=0$ for all integers $p\geq 1$, we have upon defining the operator $C_p$ by $$C_p=(-1)^{m+1}\{\sum_{j=0}^{m-1}{(-1)^{j}\left(\begin{array}{clcr}m\\j\end{array}\right)T^{p(m-j)}S^{p(m-j-1)}}\}$$ that $$C_pS^p=I,\hspace{2mm}\mbox{all integers}\hspace{2mm} p\geq1.$$ Evidently the operator $S^p$ is left invertible by $C_p$ for all integers $p\geq 1$, and $$||C_p||\leq \{1+\left(\begin{array}{clcr}m\\1\end{array}\right)+ \cdots +\left(\begin{array}{clcr}m\\m-2\end{array}\right)+\left(\begin{array}{clcr}m\\m-1\end{array}\right)\}M_1M_2<2^mM_1M_2=M$$
	for all integers $p\geq 1$. Thus, for all $x\in\X$ and integers $p\geq 1$, $$||x||=||C_pS^px||\leq M||S^px||\Longleftrightarrow (\frac{1}{M})||x||\leq ||S^px||.$$ Since already $$||S^px||\leq ||S^p||||x||\leq M_1||x||$$ for all $x\in\X$, 	
	it follows  that $S$ is similar to an invertible isometry (on an equivalent Banach space). (This is well known -- see, for example, \cite{KR}.) The proof now follows, since invertible isometries are polaroid and the polaroid property is preserved by similarities.\end{demo}
Power bounded $m$-isometric operators satisfy the property that they are isometric - see \cite[Theorem 2.4]{D} and \cite[Theorem 2.4]{DK}. Hence:
\begin{cor}\label{cor10} Power bounded $m$-isometric Banach space operators, i.e. power bounded operators $S\in\b$ such that $\triangle^m_{S^*,S}(I)=0$, are polaroid.\end{cor}
The Power bounded hypothesis on $S$ may be dropped in the case in which $\triangle^2_{S^*,S}(I)=0$ (i.e., the operator $S$ is $2$-isometric), for the reason
that invertible $2$-isometries are isometries: {\em $2$-isometric Banach space operators are polaroid.} Corollary \ref{cor10} extends to operators $S\in\B$ satisfying $\triangle^m_{S^*,CSC}(I)$ for some conjugation $C$ (i.e., to  $(m,C)$-isometries $S\in\B$). Observe that if $S$ is power bounded, then so is $CSC$ and $\sigma_a(CSC)=\overline{\sigma_a(S)}$ ($=$ complex conjugate of $\sigma_a(S)$) for every conjugation $C$. Hence:
\begin{cor}\label{cor11} Power bounded $(m,C)$-isometries $\in\B$ are polaroid.\end{cor}

{\bf Extension to $n$-quasi left $m$-invertible operators.} Theorem \ref{thm10} extends to $n$-quasi left $m$-invertible operators $S\in\B$,
$$S^{*n}\triangle^m_{T,S}(I)S^n=S^{*n}\{\sum_{j=0}^m{(-1)^{j}\left(\begin{array}{clcr}m\\j\end{array}\right)T^{m-j}S^{m-j}}\}S^n=0,$$
such that $[S,T^*]=ST^*-T^*S=0$. Letting $S$ and $T^*$ have the upper triangular representations $(2)$, it is seen that ($\sigma(S)=\sigma(S_1)\cup\{0\}$, $S_1^n(\H)=\overline{S_1^n{\H}}$, $[S_1,T^*_1]=0$ and) $$\triangle^m_{T_1,S_1}(I_1)=\sum_{j=0}^m{(-1)^{j}\left(\begin{array}{clcr}m\\j\end{array}\right)T_1^{m-j}S_1^{m-j}}=0$$
(so that $S_1$ is left $m$-invertible by $T_1$). Recall from the previous section that $S^n=\left(\begin{array}{clcr}S^n_1 & X\\ 0 & 0\end{array}\right)$, where $S_1^n$ is left $m$-invertible by $T_1^n$.  Since $S$ and $T$ power bounded implies $S_1^n$ and $T_1^n$ are power bounded, $S^n_1$ (therefore, $S_1$) is polaroid.  Hence:
\begin{thm}
	\label{thm11} Power bounded operators $S, T\in\B$ satisfying $S^{*n}\triangle^m_{T,S}(I)S^n=0$ such that $[S,T^*]=0$ are polaroid.\end{thm} \begin{demo}  Since $S=\left(\begin{array}{clcr} S_1 & S_0\\0 & S_2\end{array}\right)$, where $S_1$ is polaroid and $S_2$ is $n$-nilpotent, the proof follows from the inequalities that $\asc(S-\lambda)\leq \asc(S_1-\lambda)+\asc(S_2-\lambda)$ and $\dsc(S-\lambda)\leq \dsc(S_1-\lambda)+\dsc(S_2-\lambda)$ for all complex $\lambda$ \cite[Exexrcise 7, Page 293]{TL}.\end{demo}

\begin{rema}\label{rema10}{\em Theorem \ref{thm11} has  an $n$-quasi $m$-isometric and an $n$-quasi $(m,C)$-isometric analogue, namely:\\ {\em Power bounded
		$n$-quasi $m$-isometric operators $S\in\B$, $S^{*n}\triangle^m_{S^*,S}(I)S^n=0$,  and  power bounded $n$-quasi $(m,C)$-isometric operators $S\in\B$, $S^{*n}\triangle_{S^*,CSC}(I)S^n=0$,  such that $C=C_1\oplus C_2$  are polaroid.}\\ In particular,  $1$-quasi $2$-isometries are polaroid \cite{MP}: This follows since operators $S$ such that $S^*\triangle^2_{S^*,S}(I)S=0$ have a representation $\left(\begin{array}{clcr}S_1&X\\0&0\end{array}\right)$, where the operator $S_1$ (satisfying $\triangle^2_{S^*_1,S_1}(I_1)=0$) is polaroid. Observe here that either $\sigma(S)={\overline{\D}}$ or $\sigma(S)\subseteq\partial{\D}\cup\{0\}$.}\end{rema}

It is easily seen that for an $m$-symmetric operator $S\in\B$, $\delta^m_{S^*,CSC}(I)=0$,  $\sigma_a(S)=\overline{\sigma_a(CSC)}$ and $\lambda\in\sigma_a(S)\Longrightarrow\lambda\in\sigma_a(CS^*C)=\sigma_{su}(S)$. (Recall: $\sigma_{su}(S)=$ the surjectivity spectrum of $S$.) Hence $\sigma(S)=\sigma_a(S)\cup\sigma_{su}(S)\subseteq\sigma_{su}(S)\subseteq\sigma(S)$, i.e., $\overline{\sigma(S)}=\sigma(CSC)=\overline{\sigma_a(S)}=\overline{\sigma_{su}(S)}$. The argument of the proof of Theorem \ref{thm11} implies that {\em if the left invertible operator $S\in\B$ is $n$-quasi $m$-symmetric, $S^{*n}\delta^m_{S^*,CSC}(I)S^n=0$,  and $C=C_1\oplus C_2$, then $S$ is power bounded implies that if $S_1^n$ is polaroid, then  ($S^n$, therefore) $S$ is polaroid}.

\

For $m$-selfadjoint operators $S\in\B$, $\delta^m_{S^*,S}(I)=0$,  it is seen that  if $\lambda$ is an eigenvalue of $S$ with an eigenvector $x$ and $\overline{\mu}$ is an eigenvalue of $S^*$ with an eigenvector $y$, then $(\lambda-{\overline{\mu}})xy=0$. Hence the eigenvalues of an $m$-selfadjoint operator are real. Since $\lambda$ is a pole of $S$ implies $\lambda$ is an eigenvalue of $S$, the poles of $S$ are all real. Consider now a left invertible $n$-quasi $m$-selfadjoint operator $S\in\B$, $S^{*n}\delta^m_{S^*,S}(I)S^n=0$. Then, follow an argument similar to that above, {\em $S$ is polaroid if  the left invertible $m$-selfadjoint operator $S^n_1$ is polaroid, and this happens if and only if the isolated points of the intersection of $\sigma(S_1)$ with the real line consists of the poles of $S_1$}.

\

{\bf Selfadjoint Riesz Idempotents.} Restricting ourselves to operator $S,T\in\B$ for which $S^{*n}\triangle^m_{T,S}(I)S^n=0$ (i.e., $n$-quasi left $m$-invertible operators in $\B$)  for which $[S,T^*]=0$, in the following we consider conditions guaranteeing the self-adjointness of the Riesz idempotents $P_{\lambda}$ attached with the poles $\lambda\in\iso\sigma(S)$
of $S$. It is clear from the above that if a point $\lambda\neq 0$ is a pole of $S$, then $S$ has a matrix representation
$$S=\left(\begin{array}{clcr}
\lambda & X_1 & Y_1\\0 & S_{11} & Y_2\\0 & 0 & S_2
\end{array}\right)$$ with respect to the decomposition $\H= (S_1-\lambda)^{-1}(0)\oplus (S_1-\lambda)(\H)\oplus{S^*}^{-n}(0)$. If $x=(x_1,x_2,x_3)\in (S-\lambda)^{-1}(0)$, then (necessarily) $x_3=x_2=0$. Hence $x\in (S-\lambda)^{-1}(0)$ if and only if $x=
(x_1,0,0)$. Consider now ${(S-\lambda)^*}^{-1}(0)$. Since $(S-\lambda)^{-1}(0)\subseteq {(S-\lambda)^*}^{-1}(0)$ if and only if
$X_1^*x_1=0=Y_1^*x_1$,
$$(S-\lambda)^{-1}(0)\subseteq {(S-\lambda)^*}^{-1}(0)\Longleftrightarrow (S-\lambda)^*(S-\lambda)^{-1}(0)\subseteq\{0\}.$$
Evidently, if $(S-\lambda)^*(S-\lambda)^{-1}(0)\subset\{0\}$,
then $(S-\lambda)^{-1}(0)\subseteq {(S-\lambda)^*}^{-1}(0)$. The point $\lambda$ being
a simple pole of $S$, if $ {(S-\lambda)^*}(S-\lambda)^{-1}(0)\subseteq\{0\}$,
then ($(S-\lambda)(\H)$ is closed and)
\begin{eqnarray*}
	\H &=& (S-\lambda)^{-1}(0)\oplus (S-\lambda)(\H)=P_{\lambda}{\H}\oplus (I-P_{\lambda}){\H}\\
	&=& {(S-\lambda)^*}^{-1}(0)\oplus{(S-\lambda)^*}^{-1}(0)^{\perp}\\
	&=& (S-\lambda)^{-1}(0)\oplus (S-\lambda)^{-1}(0)^{\perp}= P_{\lambda}{\H}\oplus P_{\lambda}{\H}^{\perp}\\
	&\Longrightarrow& {P_{\lambda}\H}^{\perp}=P_{\lambda}^{-1}\H=(I-P_{\lambda})\H,
\end{eqnarray*} i.e., $P_{\lambda}$ is selfadjoint.

\

Consider now the case in which $\lambda=0$ is a pole of $S$. Then $P_{\lambda}{\H}=S^{-n}(0)$ and $S^n$ has a triangulation

$$S^n=\left(\begin{array}{clcr}S_1^n & X\\0 & 0\end{array}\right)\left(\begin{array}{clcr}{\overline{S^n{\H}}}\\{S^*}^{-n}(0)\end{array}\right),$$
where $S_1$ is invertible (since $0\in\iso\sigma(S^n)$ implies $0\notin \sigma(S^n_1)$). Since $x=(x_1,x_2)\in S^{-n}(0)$ if and only if  $x=(-S_1^{-n}Xx_2,x_2)$, $S^{-n}(0)\subseteq
{S^*}^{-n}(0)$ if and only if $Xx_2=0$, i.e., if and only if $S^n({S^*}^{-n}(0))=\{0\}$ (and then $S^{-n}(0)={S^*}^{-n}(0)$). Arguing as above, it now follows that the projection $P_0$ is selfadjoint if and only if $S^n:{S^*}^{-n}(0)\longrightarrow\{0\}$. We have proved:
\begin{pro}\label{pro10} Given an $n$-quasi left $m$-invertible operator $S\in\B$ such that $[S,T^*]=0$, the Riesz projection $P_{\lambda}$
	corresponding to a pole $\lambda\neq 0$ (resp., $\lambda=0$) of $S$ is selfadjoint if and only if $(S-\lambda)^*:(S-\lambda)^{-1}(0)\longrightarrow\{0\}$ (resp.,
	$S^n:{S^*}^{-n}(0)\longrightarrow\{0\}$).\end{pro}

\begin{rema} {\em It is immediate from the above that if $S\in\B$ is a $1$-quasi $2$-isometry, then the Riesz projection $P_{\lambda}$ corresponding to a pole $\lambda\neq 0$ (resp., $\lambda=0$)
	is selfadjoint if and only if $(S-\lambda)^*:(S-\lambda)^{-1}(0)\rightarrow\{0\}$ (resp., $S:{S^*}^{-1}(0)\rightarrow\{0\}$) ; {\em cf.} \cite
	[Theorems 2.7 and 2,8]{MP}.}\end{rema}

%%%%%%%%%%%%%%%%%%%%%%%%%%%%%%%%%%%%%%%%%%%%%%%%%%%%%%%%%
\section {\sfstp  Products }  Let $S_i, T_i\in\b$, $i=1,2$, be such that $[S_1,S_2]=0=[T_1,T_2]$ and $d^{n_i}_{T_i,S_i}(I)=0$. Then 
\begin{eqnarray*} \triangle^n_{T_1T_2,S_1S_2} &=& (L_{T_1}L_{T_2}R_{S_1}R_{S_2}-I)^n\\ &=&\{L_{T_1}(L_{T_2}R_{S_2}-I)R_{S_1}+(L_{T_1}R_{S_1}-I)\}^n\\ &=& \sum_{j=0}^n{\left(\begin{array}{clcr}n\\j\end{array}\right)\triangle_{T_2,S_2}^{n-j}L_{T_1}^{n-j}R^{n-j}_{S_1}}\triangle^j_{T_1,S_1}\end{eqnarray*}
implies 
\begin{eqnarray*} \triangle^n_{T_1T_2,S_1S_2}(I)=\sum_{j=0}^n{\left(\begin{array}{clcr}n\\j\end{array}\right)T^{n-j}_1\triangle^{n-j}_{T_2,S_2}(I)S_1^{n-j}\triangle^j_{T_1,S_1}(I)}\end{eqnarray*} 
	and \begin{eqnarray*}  \delta^n_{T_1T_2,S_1S_2} &=& (L_{T_1}L_{T_2}-R_{S_1}R_{S_2})^n\\ &=& \{L_{T_2}(L_{T_1}-R_{S_1})+(L_{T_2}-R_{S_2})R_{S_1}\}^n\\ &=& \sum_{j=0}^n{\left(\begin{array}{clcr}n\\j\end{array}\right)L^{n-j}_{T_2}\delta^{n-j}_{T_1,S_1}\delta^j_{T_2,S_2}R^j_{S_1}}	
	\end{eqnarray*}
implies  $$\delta^n_{T_1T_2,S_1S_2}(I)= \sum_{j=0}^n{\left(\begin{array}{clcr}n\\j\end{array}\right)T^{n-j}_2\delta_{T_1,S_1}^{n-j}(I)\delta^j_{T_2,S_2}(I)S^j_1}.$$
	Letting $n=m_1+m_2-1$, since $d^j_{T_2,S_2}(I)=0$ for all $j\geq m_2$ and $d^{m_1+m_2-1-j}_{T_1S_1}(I)=0$ for all $j\leq m_2-1$ (implies $m_1+m_2-1-j\geq m$) , we have:
	\begin{lem}
		\label{lem300} If $S_i, T_i\in\b$, $i=1,2$, are such that  $[S_1,S_2]=0=[T_1,T_2]$ and $d^{m_i}_{T_i,S_i}(I)=0$, then $d^{m_1+m_2-1}_{T_1T_2,S_1S_2}(I)=0$. \end{lem}

 The following
theorem is an $n(S)$-quasi $[m,d]$-version of these results. (Recall here that the  operators {\em $S_1, T_1\in\b$ are $n(S)$-quasi $[m,d]$-intertwined for an operator $S\in\B$ if  $S^{*n}d^m_{T_1,S_1}(I)S^n=0$}.)
\begin{thm}
	\label{thm01} If $S^{*n}d^{m_i}_{T_i,S_i}(I)S^n=0$, $i=1,2$, for some operators $S, S_1, S_2, T_1, T_2\in\B$ such that $[S,S_i]=0=[S,T^*_i]$  and $[S_1,S_2]=0=[T_1,T_2]$, then $S^{*n}d^{m_1+m_2-1}_{T_1T_2,S_1S_2}(I)S^n=0$ (i.e., $T_1T_2$ and $S_1S_2$ are $n(S)$-quasi $[m_1+m_2-1, d]$-intertwined).\end{thm}
\begin{demo} The hypotheses imply that the operators $S, S_i$ and $T_i^*$ have the upper triangular matrix representations
	
	\begin{eqnarray*}
		& & S=\left(\begin{array}{clcr}S_{01} & S_{00}\\0 & S_{02}\end{array}\right), \hspace{2mm} S_i=\left(\begin{array}{clcr}S_{i1} & S_{i0}\\0 & S_{i2}\end{array}\right),\\ & & T_i^*=\left(\begin{array}{clcr}T_{i1}^* & T^*_{i0}\\0 & T^*_{i2}\end{array}\right); i=1,2,\end{eqnarray*} with respect to the decomposition $\H=\overline{S^n(\H)}\oplus {S^*}^{-n}(0)$ of $\H$. The hypothesis $S^{*n}d^{m_1}_{T_1,S_1}(I)S^n=0$ implies $d^{m_1}_{T_{11},S_{11}}(I_1)=0$ and the hypothesis $S^{*n}d^{m_2}_{T_2,S_2}(I)S^n=0$ implies $d^{m_2}_{T_{21},S_{21}}(I_1)=0$. Hence, since the hypothesis $[S_1,S_2]=0=[T_1,T_2]$ implies $[S_{11},S_{21}]=0=[T_{11},T_{21}]$,
	Lemma \ref{lem300} implies $d^{m_1+m_2-1}_{T_{11}T_{21},S_{11}S_{21}}(I_1)=0$. Finally, since
	$$\sum_{j=0}^{m_1+m_2-1}{(-1)^{j}\left(\begin{array}{clcr}m_1+m_2-1\\j\end{array}\right) (T_1T_2)^{m_1+m_2-1-j}(S_1S_2)^{m_1+m_2-1-j}}=\left(\begin{array}{clcr} 0 & Z_1\\Z_2 & Z_3\end{array}\right)$$ for some operators $Z_i$ ($i=1,2,3$), and $S^n=\left(\begin{array}{clcr}S^n_{01} & X\\0 & 0\end{array}\right)$ for some operator $X$, with respect to $\H=\overline{S^n(\H)}\oplus {S^*}^{-n}(0)$,
	$$ S^{*n}d^{m_1+m_2-1}_{T_1T_2,S_1S_2}(I)S^n=0,$$
	i.e., $T_1T_2$ and $S_1S_2$ are $n(S)$-quasi $[m_1+m_2-1,d]$-intertwined.\end{demo}

\begin{rema}\label{rem01} {\em (i) Recall that $T$ is a strict left $m$-inverse of $S$ if $\triangle^m_{S,T}(I)=0$ but $\triangle^{m-1}_{S,T}(I)\neq 0$ \cite{{DM}, {CG}}. Letting $m_1=1$ in $\triangle^{m_1}_{T_1,S_1}(I)=0$ (so that $T_1$ is a left $1$-inverse of $S_1$, i.e., $T_1S_1=I$), it follows that  $T_1T_2$ is a strict left $m_2$-inverse of $S_1S_2$ if and only if $\triangle^{m_2-1}_{T_2,S_2}(I)\neq 0$ \cite[Theorem 13]{CG}, i.e., if and only if $T_2$ is a strict left $m_2$-inverse of $S_2$. Theorem \ref{thm01} does not extend to  $n(S)$-quasi strict $[m_1+m_2-1, d]$-intertwinings. Thus, given
	$T_1$ an $n(S)$-quasi left $1$-inverse of $S_1$ (i.e., $S^{*n}\triangle_{T_1,S_1}(I)S^n=0$) and $T_2$ a strict left $m$-inverse of $S_2$ (i.e., $\triangle^{m}_{T_2,S_2}(I)=0$ and $\triangle^{m-1}_{T_2,S_2}(I)\neq 0$), $T_1T_2$ may not be an $n(S)$-quasi strict left $m$-inverse of $S_1S_2$. To see this, consider operators $S_i$ and $T_i$ satisfying the commutativity hypotheses of Theorem
	\ref{thm01} such that $T_{11}$ is left $1$-inverse of $S_{11}$, $T_{21}$ is a left $(m-1)$-inverse of $S_{21}$ and $T_{22}$ is a strict left $m$-inverse of $S_{22}$. Define $S_i$ and $T_i$ by $$S_1=S_{11}\oplus I, S_2=S_{21}\oplus S_{22}, T_1=T_{11}\oplus I,\hspace{2mm}\mbox{and}\hspace{2mm} T_2=T_{21}\oplus T_{22}$$
	(with respect to the decomposition $\H=\overline{S^n(\H)}\oplus {S^*}^{-n}(0)$ of $\H$). Then $T_1$ is an $n(S)$-quasi left $1$-inverse of $S_1$, $T_2$ is a strict left $m$-inverse of $S_2$, and $T_1T_2$ is not an $n(S)$-quasi strict left $m$-inverse of $S_1S_2$.
	
	\
	
	(ii) Trivially, one may replace $n(S)$-quasi by $n(S_iS)$-quasi, $i=1,2$, in the conclusion of Theorem \ref{thm01}.} \end{rema}

Given Hilbert spaces $\H_i$, $i=1,2$, let $H_1{\overline{\otimes}}\H_2$ denote the completion, endowed with a reasonable uniform cross-norm,
of the algebraic tensor product $\H_1\otimes\H_2$ and, for  $A_i\in B(\H_i)$, $i=1,2$, let  $A_1\otimes A_2\in B(H_1{\overline{\otimes}}\H_2)$ denote the tensor product of $A_1$ and $A_2$.  Theorem \ref{thm01} applies to tensor products of $n$-quasi left $m$-invertible, $m$-isometric and $(m,C)$-isometric operators.
Let $A_i, B_i$ ($i=1,2$) and $S,T$ be operators in $\B$.
\begin{cor}\label{cor01} If $A^{*n}_1d^{m_1}_{B_1,A_1}(I)A^n_1=0=d^{m_2}_{B_2,A_2}(I)$  and
	$[A_1,B^*_1]=0$, then $(A_1\otimes A_2)^{*n}d^{m_1+m_2-1}_{B_1\otimes B_2, A_1\otimes A_2}(I\otimes I)(A_1\otimes A_2)^n=0$. \end{cor}
\begin{demo} Define the operators $S, S_i$ and $T_i$, $i=1,2$, by
	$$S=S_1=A_1\otimes I, T_1=B_1\otimes I, S_2=I\otimes A_2\hspace{2mm}\mbox{and}\hspace{2mm} T_2=I\otimes B_2.$$ Then, since $[A_1,B^*_1]=0$,
	$$[S_1,S_2]=0=[T_1,T_2]\hspace{2mm}\mbox{and}\hspace{2mm} [S,T_i^*]=0=[S,S_i]$$
	($i=1,2$). Theorem \ref{thm01} applies to prove
	\begin{eqnarray*}	 (A_1^*\otimes I)^n d^{m_1+m_2-1}_{B_1\otimes B_2,A_1\otimes A_2}(I\otimes I)(A_1\otimes I)^n=0.\end{eqnarray*} Multiplying by $(I\otimes B_1^*)^n$ on the left and by $(I\otimes B_1)^n$ on the right, the
	proof follows.\end{demo}
Translated to $(m,C)$-isometric operators, Theorem \ref{thm10} and Corollary \ref{cor10} imply the following.
\begin{cor}\label{cor02} Given conjugations $C$ and $D$, if:
	\vskip3pt\noindent (i) $S,T\in\B$ are commuting operators such that  $S^{*n}\triangle^{m_1}_{S^*,CSC}(I)S^n=0=\triangle^{m_2}_{T^*,DTD}(I)$,  $[S,CSC]=0=[S,DTD]$ and $[T,CSC]=0=[DTD,CSC]$, then
	
	\begin{eqnarray*}  (ST)^{*n}\triangle^{m_1+m_2-1}_{S^*T^*,CSCDTD}(I)(ST)^n &= & (ST)^{*n}\{\sum_{j=0}^{m_1+m_2-1}{(-1)^{j}\left(\begin{array}{clcr}m_1+m_2-1\\j\end{array}\right) }\times\\  &\times& { (ST)^{*(m_1+m_2-1-j)} (CSCDTD)^{m_1+m_2-1-j}\}(ST)^n}\\ &=&0.\end{eqnarray*}
	In particular, if $C=D$, then $$(ST)^{*n}\triangle^{m_1+m_2-1}_{S^*T^*,CSTC}(I)(ST)^n=0$$ (i.e.,  $ST$ is $n$-quasi $(m_1+m_2-1, C)$-isometric).
	
	\vskip3pt\noindent (ii) $A^{*n}\triangle^{m_1}_{A^*,CAC}(I)A^n=0=B^{*n}\triangle^{m_2}_{B^*,DBD}(I)B^n$  and $[A,CAC]=0$, then $$(A\otimes B)^{*n}\triangle^{m_1+m_2-1}_{A^*\otimes B^*, (CAC\otimes DBD)}(I\otimes I)(A\otimes B)^n=0$$ (i.e.,  $A\otimes B$ is $n$-quasi $(m_1+m_2-1,C\otimes D)$-isometric).
\end{cor}
\begin{demo} (i) If we define $S_i$ and $T_i$, $i=1,2$, by $S_1=CSC$, $S_2=DTD$, $T_1=S^*$ and $T_2=T^*$, then $S,S_i$ and $T_i$ ($i=1,2$)
	satisfy the hypotheses of Theorem \ref{thm01}. Hence the proof of (i). The proof of (ii) is evident.\end{demo}
\

Corollary 4.4 generalizes \cite[Theorem 2.3]{ACL} (proved for the case $n=0$ and $C=D$), and Corollaries 2.1, 3.5 and Proposition 3.5 (proved for the cases $n=2,3$ of part (ii) of our Corollary 4.4)  of \cite{ACL}.

\

Corollary \ref{cor02} takes the following simpler form for $m$-isometries.
\begin{cor}\label{cor03} Given operators $S, T\in\B$ such that $S^{*n}\triangle^{m_1}_{S^*,S}(I)S^n=0=\triangle^{m_2}_{T^*,T}(I)$ (i.e.,  $S$ is $n$-quasi $m_1$-isometric and $T$ is $m_2$-isometric):
	\vskip3pt\noindent(i) if $[S,T]=0$, then $(ST)^{*n}\triangle^{m_1+m_2-1}_{S^*T^*,ST}(I)(ST)^n=0$ (i.e., $ST$ is $n$-quasi $(m_1+m_2-1)$-isometric);
	
	\vskip3pt\noindent (ii) $(S\otimes T)^{*n}\triangle^{m_1+m_2-1}_{S^*\otimes T^*,S\otimes T}(I\otimes I)(S\otimes T)^n=0$  (i.e., $S\otimes T$ is $n$-quasi $(m_1+m_2-1)$-isometric).
\end{cor}

\

 A version of Corollary \ref{cor03} holds for $m$-selfadjoint and $m$-symmetric operators.

\begin{cor}\label{cor04} Let $S,T\in\B$ satisfy $[S,T]=0$ and let $C$ be a conjugation of $\H$. If:

	\vskip4pt\noindent (i) $S$ is $n$-quasi $m_1$-selfadjoint and $T$ is $m_2$ selfadjoint, then $ST$ is
	$n$-quasi $(m_1+m_2-1)$-selfadjoint (i.e., $S^{*n}\delta^{m_1}_{S^*,S}(I)S^n)=0=\delta^{m_2}_{T^*,T}(I)=0$ implies $(ST)^{*n}\delta^{m_1+m_2-1}_{S^*T^*,ST}(I)(ST)^n=0$);

\vskip4pt\noindent (ii)  $S$ is $n$-quasi $m_1$-symmetric with the symmetry implemented by the conjugation $C$, $T$ is $m_2$-symmetric with the symmetry implemented by the conjugation $C$ and $[S,CSC]=0$, then $ST$ is $n$-quasi $m_1+m_2-1$-symmetric with the symmetry implemented by the conjugation $C$ (i.e.,
$S^{*n}\delta^{m_1}_{S^*,CMC}(I)S^n=0=\delta^{m_2}_{T^*,CTC}(I)=0$ and $[S,CSC]=0$ implies $(ST)^{*n}\delta^{m_1+m_2-1}_{S^*T^*,CSTC}(I)(ST)^n$);
	
	\vskip4pt\noindent (iii)  $S$ is $n$-quasi $m_1$-selfadjoint and  $T$ is $m_2$-selfadjoint, then $S\otimes T$ is $n$-quasi $(m_1+m_2-1)$-selfadjoint;

\vskip4pt\noindent (iv)  $S$ is $n$-quasi $m_1$-symmetric and $T$ is $m_2$-symmetric (with the symmetry implemented by the conjugation $C$ for $S$ and $T$), then $S\otimes T$ is $n$-quasi $(m_1+m_2-1)$-symmetric (with the symmetry implemented by the conjugation $C$).
\end{cor}

\section {\sfstp Perturbation by Nilpotents. }  Gu \cite[Theorem 2]{CG} proves that if $T\in\b$ is a left (right) $m$-inverse of $S\in\b$
and $N\in\b$ is an $n$-nilpotent which commutes with $T$, then $T+N$ is a left (resp., right)  $(m+n-1)$-inverse of $S$. Consequently, {\em
	If $T$ is a left $m$-inverse of $S$, $N_1$ is an $n_1$-nilpotent which commutes with $T$ and $N_2$ is an $n_2$-nilpotent which commutes with $S$, then $T+N_1$ is a left $(m+n_1+n_2-2)$-inverse of $S+N_2$.} Translated to $m$-isometric (and $(m,C)$-isometric) operators $S$, this implies: {\em If $N\in\B$ is an n-nilpotent operator which commutes with $S$, then $S+N$ is an $(m+2n-2)$-isometric \cite{BMMN} (resp., $(m+2n-2,C)$-isometric \cite{CLM}) operator.} A similar result holds for $m$-selfadjoint and $(m,C)$-symmetric operators \cite{{CKL}, {L}}. In the following we consider perturbation by commuting nilpotents of operators $S,T\in \b$ satisfying $d^m_{T,S}(I)=0$, and using an elementary argument we prove:
\begin{thm}\label{thm301} If $d^m_{T,S}(I)=0$ and $N\in\b$ is an $n$-nilpotent operator satisfying $[S,N]=0$, then $d^{m+n-1}_{T,S+N}(I)=0$.\end{thm}
\begin{demo} We start by proving that \begin{eqnarray*} & & \triangle^p_{T,S+N}(I)=\sum_{j=0}^p{\left(\begin{array}{clcr}p\\p-j\end{array}\right)T^j\triangle^{p-j}_{T,S}(I)N^j,}\hspace{2mm}\mbox{and}\\ & &   \delta^p_{T,S+N}(I)=\sum_{j=0}^p{(-1)^j\left(\begin{array}{clcr}p\\j\end{array}\right)\delta^{p-j}_{T,S}(I)N^j}.\end{eqnarray*} The proof is by induction. Both the equalities being true for $p=1$, assume their validity for some $k>1$. Then
\begin{eqnarray*} & & \triangle^{k+1}_{T,S+N}(I)= \triangle_{T,S}( \triangle^{k}_{T,S+N}(I))+T \triangle^{k}_{T,S+N}(I)N\\ &=&   \triangle^{k+1}_{T,S}(I)+\{\left(\begin{array}{clcr}k\\k\end{array}\right)+\left(\begin{array}{clcr}k\\k-1\end{array}\right)\}T\triangle^k_{T,S}N+\{\left(\begin{array}{clcr}k\\k-1\end{array}\right)+\left(\begin{array}{clcr}k\\k-2\end{array}\right)\}T^2\triangle^{k-1}_{T,S}N^2\\   & & \cdots+ \{\left(\begin{array}{clcr}k\\1\end{array}\right)+\left(\begin{array}{clcr}k\\0\end{array}\right)\}T^k\triangle_{T,S}N^k+\left(\begin{array}{clcr}k\\0\end{array}\right)T^{k+1}N^{k+1}\\ &=& \sum^{k+1}_{j=0}{\left(\begin{array}{clcr}k+1\\j\end{array}\right) T^j\triangle^{k+1-j}_{T,S}(I)N^j}, \hspace{2mm}\mbox{and}\end{eqnarray*}
\begin{eqnarray*} & & \delta^{k+1}_{T,S+N}(I)= \delta_{T,S}(\delta^k_{T,S+N}(I))-\delta^k_{T,S+N}(I)N\\ &=&  \delta^{k+1}_{T,S}(I)+\{(-1)\left(\begin{array}{clcr}k\\0\end{array}\right)-\left(\begin{array}{clcr}k\\1\end{array}\right)\}\delta^k_{T,S}N+\{(-1)^2\left(\begin{array}{clcr}k\\2\end{array}\right)-(-1)\left(\begin{array}{clcr}k\\1\end{array}\right)\}\delta^{k-1}_{T,S}N^2\\   & & \cdots+ \{(-1)^k\left(\begin{array}{clcr}k\\k\end{array}\right)-(-1)^{k-1}\left(\begin{array}{clcr}k\\k-1\end{array}\right)\}\delta_{T,S}N^k-(-1)^k\left(\begin{array}{clcr}k\\k\end{array}\right)N^{k+1}\\ &=& \sum^{k+1}_{j=0}{(-1)^j\left(\begin{array}{clcr}k+1\\j\end{array}\right) \delta^{k+1-j}_{T,S}(I)N^j}.\end{eqnarray*} Recall now that $d^m_{T,S}(I)=0$ implies $d^t_{T,S}(I)=0$ for all integers $t\geq m$. Hence, since $N^j=0$ for all $j\geq n$, $d^p_{T,S+N}(I)=0$ for all $p$ such that $p-n+1\geq m$ (in particular, if $p=m+n-1$).\end{demo}
Trivially, $d^m_{T,S}(I)=0$ if and only if $d^m_{S^*,T^*}(I)=0$ (where we have used $I$ to denote the identity of both $\b$ and $B(\X^*)$). Hence:

 \begin{cor}\label{cor301}  If $d^m_{T,S}(I)=0$ and $N_i\in\b$ ($i=1,2$) are $n_i$-nilpotent operators satisfying $[S,N_1]=0=[T,N_2]$, then $d^{m+n_1+n_2-2}_{T+N_2,S+N_1}(I)=0$.\end{cor}

For perturbation by commuting nilpotents of $n$-quasi $[m, d]$-operators (i.e., operators $S,T\in\B$ such that $S^{*n}d^m_{T,S}(I)S^n=0$), we have the following.
\begin{thm}\label{thm30}  Suppose that $S^{*n}d^m_{T,S}(I)S^n=0$ for some operators $S,T\in\B$ and integers $m,n\geq 1$.  If  $N_i\in\B$, $i=1,2$, are $n_i$-nilpotent operators such that $[S,N_1]=0=[S,T^*]$ and $[N_2,T]=0=[N_2^*,S]$, then  $$(S^*+N_1^*)^{n+n_1-1}d^{m+n_1+n_2-2}_{T+N_2,S+N_1}(I)(S+N_1)^{n+n_1-1}=0.$$\end{thm}
\begin{demo} Letting $S$ and $T^*$ have the upper triangular representations $(3)$ of Section 2, it follows from the hypotheses that
	$N_1$ and $N_2$ have the upper triangular representations
	$$N_1=\left(\begin{array}{clcr}N_{11} & N_{10}\\0 & N_{12}\end{array}\right)\hspace{2mm}\mbox{and}\hspace{2mm} N^*_2=\left(\begin{array}{clcr} N^*_{21} & N^*_{20}\\0 & N^*_{22}\end{array}\right)$$
	(with respect to the decomposition $\H=\overline{S^n(\H)}\oplus {S^*}^{-n}(0)$), where $$N_{11}^{n_1}=N_{12}^{n_1}=0=N_{21}^{n_2}=N_{22}^{n_2}\hspace{2mm}\mbox{and}\hspace{2mm} [N_{11},S_1]=0=[N_{21},T_1].$$
	The hypothesis  $S^{*n}d^m_{T,S}(I)S^n=0$  implies  $d^m_{T_1,S_1}(I_1)=0$. Hence ,
	
	$$ d^{m+n_1+n_2-2}_{T_1+N_{21},S_1+N_{11}}(I_1)=0.$$
	This, since
	 $$(S+N_1)^{n+n_1-1}=\left(\begin{array}{clcr}(S_1+N_{11})^{n+n_1-1} & Z\\0 & 0\end{array}\right)$$ (for some operator $Z$) and
$$d^{m+n_1+n_2-2}_{T+N_2,S+N_1}(I)=\left(\begin{array}{clcr} 0&Z_1\\Z_2&Z_3\end{array}\right)$$	
for some operators $Z_i$ ($i=1,2,3$), implies
	\begin{eqnarray*}
		(S^*+N_1^*)^{n+n_1-1}d^{m+n_1+n_2-2}_{T+N_2,S+N_1}(I)(S+N_1)^{n+n_1-1}=0.\end{eqnarray*} This completes the proof.	
\end{demo}

More can be said in the case in which $T=S^*$ (i.e., when $S$ is $n$-quasi $m$-isometric \cite{ACL}).
\begin{cor}\label{cor30} Given an operator $S\in\B$ such that $S^{*n}\triangle^m_{S^*,S}(I)S^n=0$, let $N\in\B$ be an $n_1$-nilpotent operator such that $[S,N]=0$. Then:
	\vskip3pt\noindent (i) $S^{*(n+n_1-1)}\triangle^{m+2n_1-2}_{S^*+N^*,S+N}(I)(S+N)^{n+n_1-1}=0$  (i.e., $(S+N)$ is an $(n+n_1-1)$-quasi $(m+2n_1-2)$-isometric operator).
	\vskip3pt\noindent(ii) If $S_1=S|_{\overline{S^n(\H)}}$ has a dense range (or, $S_1^*$ has SVEP at $0$), then $(S+N)^{n+n_1-1}$ is similar to the operator $(S_1+N_1)^{n+n_1-1}\oplus 0$.
	\vskip3pt\noindent (iii) There exists a positive operator $Q$ and an operator $A$ similar to $S^{n+n_1-1}$ such that $\triangle^{m+2n_1-2}_{A^*,A}(Q)=0$ (i.e., $A$ is $(m+2n_1-2,Q)$-isometric. Furthermore, if also $S$ is left invertible, then $S^{n+n_1-1}$ is similar to an $(m+2n_1-2)$-isometric operator.                                                                .\end{cor}
\begin{demo} The proof of (i) follows from Theorem \ref{thm30}. To prove (ii), we start by observing that if we let $N=
	\left(\begin{array}{clcr}N_1 & N_0\\0 & N_2\end{array}\right)$ (with respect to the decomposition $\H=\overline{S^n(\H)}\oplus {S^*}^{-n}(0)$), then ($N_1^n=N_2^n=0$ and $S_2^n=0$ in the corresponding representation $(3)$  for $S$)
	$$(S+N)^{n+n_1-1}=\left(\begin{array}{clcr}(S_1+N_1)^{n+n_1-1} & X\\0 & 0\end{array}\right)$$
	for some operator $X$. The operators $S_1$ and $N_1$ commute, and $S_1$ is left invertible (since $S_1$ is left $m$-invertible). Hence, since $\sigma_a(S_1+N_1)\subseteq
	\sigma_a(S_1)+\sigma_a(N_1)=\sigma_a(S_1)$, $S_1+N_1$ is left invertible. Define the operator $E\in\B$ by $E=\left(\begin{array}{clcr}
	(S_1+N_1)^{n+n_1-1} & X\\0 & 1\end{array}\right)$; then ( since either of the hypotheses $S_1$ has a dense range and $S^*_1$ has SVEP at $0$ implies) $E$ is invertible with $$E^{-1}= \left(\begin{array}{clcr}(S_1+N_1)^{-(n+n_1-1)} &
	-(S_1+N_1)^{-(n+n_1-1)}X\\0 & 1\end{array}\right).$$ If we now define $A\in \B$ by $A=(S_1+N_1)^{n+n_1-1}\oplus 0$, then
	$(S+N)^{n+n_1-1}=E^{-1}AE$. To prove (iii),  we start by observing from the proof of Theorem \ref{thm30} that the current hypotheses imply  $(S_1+N_1)^p$ is $(m+2n_1-2)$-isometric and $(S+N)^p$ is $(n+n_1-1)$-quasi $(m+2n_1-2)$-isometric for all integers $p\geq 1$. Choose $p=n+n_1-1$ and let $(S_1+N_1)^{n+n_1-1}$ have the polar decomposition $(S_1+N_1)^{n+n_1-1}=U_1P_1$ (so that $U_1$ is an isometry and $P_1$ is positive invertible). Let $(S+N)^{n+n_1-1}=\left(\begin{array}{clcr}U_1P_1&X\\0&0\end{array}\right)$ and argue as in the proof of Proposition \ref{pro00}. Then, upon defining $Q\geq 0$ as in the proof of Proposition \ref{pro00} and letting $m+2n_1-2-j=t$,
	\begin{eqnarray*} & & (S+N)^{*(n+n_1-1)}\{\sum_{j=0}^{m+2n_1-2}{(-1)^{j}\left(\begin{array}{clcr}m+2n_1-2\\j\end{array}\right)(S+N)^{*t(n+n_1-1)}}\times\\ &\times& {(S+N)^{t(n+n_1-1)}}\}(S+N)^{n+n_1-1}=0\\ &\Longleftrightarrow& \sum_{j=0}^{m+2n_1-2}{(-1)^{j}\left(\begin{array}{clcr}m+2n_1-2\\j\end{array}\right)\left(\begin{array}{clcr}U^*_1P_1&0\\X^*P_1&0\end{array}\right)^t Q \left(\begin{array}{clcr}P_1U_1&P_1X\\0&0\end{array}\right)^t}=0.\end{eqnarray*}
Now define the operator $A$ by $A=\left(\begin{array}{clcr}P_1U_1&P_1X\\0&0\end{array}\right)$. Then $A$ is $(m+2n_1-2,Q)$-isometric and $S^{n+n_1-1}=P^{-1}AP$, where $P=P_1\oplus I_2$. To complete the proof, assume now that $S$ is left invertible. Then $P$ and $Q$ are invertible positive operators, $B=Q^{\frac{1}{2}}AQ^{\frac{-1}{2}}$ is $(m+2n_1-2)$-isometric and $S^{n+n_1-1}=E^{-1}BE$, $E=Q^{\frac{1}{2}}P$.\end{demo}

The corresponding result for $n$-quasi $(m,C)$-isometries $S$, $S^{*n}\triangle^m_{S^*,CSC}(I)S^n=0$, such that $C=C_1\oplus C_2$ (with respect to the decomposition $\H=\overline{S^n(\H)}\oplus {S^*}^{-n}(0)$) is the following.   Define the operator $M$ (as before) by $M=\left(\begin{array}{clcr}U_1&X\\0&0\end{array}\right)$, where the isometry $U_1$ and the operator $X$ are as in the polar decomposition (above) of $S^{n+n_1-1}$.
\begin{cor}\label{cor31} Let $S\in\B$ be an  $n$-quasi $[m,C]$-isometry such that $C=C_1\oplus C_2$ with respect to the decomposition $\H=\overline{S^n(\H)}\oplus {S^*}^{-n}(0)$. If $N\in\B$ is an $n_1$-nilpotent operator such that $[S,N]=0$, then:
	\vskip4pt\noindent (i) $S+N$ is $(n+n_1-1)$-quasi $(m+2n-1,C)$-isometric.

	\vskip4pt\noindent(ii) $(S+N)^{n+n_1-1}$ is similar to $(S_1+N_1)^{n+n_1-1}\oplus 0$, $S_1=S|_{\overline{S^n(\H)}}$ and $N_1=N|_{\overline{S^n(\H)}}$, whenever $S_1$ has a dense range (or $S^*_1$ has SVEP at $0$).

	\vskip4pt\noindent (iii) If also $[C,M]=0$, then $(S+N)^{n+n_1-1}$ is similar to an $(m+2n_1-2,C)$-isometry.\end{cor}
\begin{demo} The hypothesis \begin{eqnarray*} & & S^{*n}\triangle^m_{S^*,CSC}(I)S^n=0\Longrightarrow \triangle^m_{S^*_1,C_1S_1C_1}(I_1)=0\\ &\Longleftrightarrow& \triangle^m_{C_1S_1^*C_1,S_1}(I_1)=0 \Longrightarrow \triangle^{m+n_1-1}_{C_1S^*_1C_1,S_1+N_1}(I_1)=0\\ &\Longleftrightarrow& \triangle^{m+n-1}_{C_1(S_1+N_1)^*C_1,S_1}(I_1)=0\Longrightarrow \triangle^{m+2n_1-2}_{C_1(S_1+N_1)^*C_1,S_1+N_1}(I_1)=0\\ &\Longleftrightarrow&  \triangle^{m+2n_1-2}_{S_1^*+N_1^*,C_1(S_1+N_1)C_1}(I_1)=0\\ &\Longrightarrow& (S+N)^{*(n+n_1-1)}\triangle^{m+2n_2-2}_{S^*+N^*,C(S+N)C}(I)S^{n+n_1-1}=0.\end{eqnarray*} This proves (i). The proof of (ii) follows from the proof of Corollary \ref{cor30}, and the proof of (iii) follows from the argument of the proof of Corollary \ref{cor30} and Proposition \ref{pro110} applied to
\begin{eqnarray*} & &  (S+N)^{*(n+n_1-1)}\triangle^{m+2n_2-2}_{S^*+N^*,C(S+N)C}(I)S^{n+n_1-1}=0 \\ &\Longrightarrow&  \sum_{j=0}^{m+2n_1-2}{(-1)^{j}\left(\begin{array}{clcr}m+2n_1-2\\j\end{array}\right) \left(\begin{array}{clcr}U^*_1P_1&0\\X^*P_1&0\end{array}\right)^t} \times\\
&\times& {\left(\begin{array}{clcr}I_1 & U^*_1C_1XC_2\\X^*U_1C_1U_1 & X^*C_1XC_2\end{array}\right) \left(\begin{array}{clcr}C_1P_1U_1C_1&C_1P_1XC_2\\0&0\end{array}\right)^t}=0,\end{eqnarray*} where $t=m+2n_1-2-j$.This completes the proof.\end{demo}

We remark in closing that Corollaries \ref{cor30} and \ref{cor31} have an $m$-selfadjoint and $m$-symmetric operators version. For example, if $S\in\B$ is satisfies $S^{*n}\delta^m_{S^*,S}(I)S^n=0$ and $N\in\B$ is an $n_1$-nilpotent which commutes with $S$, then:

\vskip4pt\noindent (i) $S_1+N_1$, where $S_1=S|_{\overline{S^n(\H)}}$ and $N_1=N|_{\overline{S^n(\H)}}$, satisfies $\delta^{m+2n_1-2}_{S^*_1+N^*_1,S_1+N_1}(I_1)=0$;

\vskip4pt\noindent (ii) $(S^*+N^*)^{n+n_1-1}\delta^{m+2n_1-2}_{S^*+N^*,S+N}(I)(S+N)^{n+n_1-1}=0$;

\vskip4pt\noindent  (iii) if also $S$ is left invertible, then $(S+N)^{n+n_1-1}$ is similar to an $(m+2n_1-2)$-selfadjoint operator.

\vskip3pt\noindent We leave the proof of the above, and the formulation of the corresponding result for $m$-symmetric operators (for which $C=C_1\oplus C_2: {\overline{S^n(\H)}}\oplus {S^*}^{-n}(0)\longrightarrow {\overline{S^n(\H)}}\oplus {S^*}^{-n}(0))$) to the reader.

%%%%%%%%%%%%%%%%%%%%%%%%%%%%%%%%%%%%%%%%%%%%%%%%%%%%%%%%%  REFERENCES

\vskip10pt \noindent\normalsize\rm B.P. Duggal, 8 Redwood Grove, London W5 4SZ, England (U.K.).\\
\noindent\normalsize \tt e-mail: bpduggal@yahoo.co.uk

\vskip6pt\noindent \noindent\normalsize\rm I. H. Kim, Department of
Mathematics, Incheon National University, Incheon, 406-772, Korea.\\
\noindent\normalsize \tt e-mail: ihkim@inu.ac.kr

\

(Note from author 1. This is the corrected version of the first author's retracted paper "On $n$-quasi left $m$-invertible operators", arXiv:1812.00221[math.FA] and  Functional Analaysis  Approximation and Computation 11(1)(2019),21-37.)

\end{document}